\documentstyle[12pt,amssymb]{article} 

\font\csc=cmcsc10 scaled\magstep1
                
\font\teneufm=eufm10
\font\seveneufm=eufm7
\font\fiveeufm=eufm5
\newfam\eufmfam
\textfont\eufmfam=\teneufm
\scriptfont\eufmfam=\seveneufm
\scriptscriptfont\eufmfam=\fiveeufm

\let\goth\frak
\newcommand{\slth}{\widehat{\goth{sl}}_2}
\newcommand{\slt}{\goth{sl}_2}
\newcommand{\slnh}{\widehat{\goth{sl}}_n}
\newcommand{\sln}{\goth{sl}_n}

\newcommand{\Aqp}[1]{{{\cal A}_{q,p^{#1}}}}

\newcommand{\Bqla}{{{\cal B}_{q,\lambda}}}

\newcommand{\sct}{{\scriptstyle}}

\font\fourteeneufm=eufm10 scaled\magstep2    
   
\newcommand{\slthbig}{\widehat{\mbox{\fourteeneufm sl}}_2}  
\newcommand{\slnhbig}{\widehat{\mbox{\fourteeneufm sl}}_n}  
\newcommand{\sltbig}{\mbox{\fourteeneufm sl}_2}  
\newcommand{\gbig}{\mbox{\fourteeneufm g}} 

\newcommand{\Z}{{\Bbb Z}} 
\newcommand{\C}{{\Bbb C}} 

\newcommand{\B}{{\cal B}}
\newcommand{\cR}{{\cal R}}
\newcommand{\ctR}{\widetilde{\cal R}}

\newcommand{\la}{\lambda}
\newcommand{\La}{\Lambda}
\newcommand{\ve}{\varepsilon}

\newcommand{\De}{{\Delta}}

\newcommand{\ze}{{\zeta}}

\newcommand{\cC}{{\cal C}}

\newcommand{\bR}{{\overline{R}}}
\newcommand{\Rck}{{\check{R}}}

\newcommand{\nn}{\nonumber}
\newcommand{\eqref}[1]{(\ref{#1})}
\newcommand{\bea}{\begin{eqnarray}}
\newcommand{\ena}{\end{eqnarray}}
\newcommand{\be}{\begin{eqnarray*}}
\newcommand{\en}{\end{eqnarray*}}
\newcommand{\lb}[1]{\label{#1}}

\newcommand{\smod}{~\hbox{\scriptsize mod}~}

\newcommand{\kuru}{\curvearrowleft} 
\newcommand{\id}{\mbox{\rm id}}
\newcommand{\Ad}{\mbox{Ad}}

\newcommand{\BW}[5]
{\left(\begin{array}{cc}#1 & #2 \cr #3 & #4 \cr\end{array}
\Biggl| #5\right)}
\newcommand{\hg}[5]
{{}_2\phi_1\left({{#1\,\,\,#2}\atop #3};#4,#5\right)}
\newcommand{\Psis}{{\Psi^*}}

\newcommand{\qed}{\hfill \fbox{}\medskip}
\newcommand{\proof}{\medskip\noindent{\it Proof.}\quad }

\newtheorem{thm}{Theorem}[section]
\newtheorem{prop}[thm]{Proposition}
\newtheorem{lem}[thm]{Lemma}
\newtheorem{dfn}[thm]{Definition}
\newcommand{\ignore}[1]{}

\newcommand{\maprightu}[2]
  {\smash{\mathop{\hbox to #1{\rightarrowfill}}\limits^{#2}}}
\newcommand{\maprightd}[2]
  {\smash{\mathop{\hbox to #1{\rightarrowfill}}\limits_{#2}}}
\newcommand{\mapdownl}[1]
  {\Bigg\downarrow\llap{$\vcenter{\hbox{$\scriptstyle#1\;\;\;$}}$}}
\newcommand{\mapdownr}[1]
  {\Bigg\downarrow\rlap{$\vcenter{\hbox{$\scriptstyle#1$}}$}}

\begin{document}
\renewcommand{\thefootnote}{\fnsymbol{footnote}}
\font\csc=cmcsc10 scaled\magstep1

{\baselineskip=14pt
 \rightline{
 \vbox{
       December 1997 \hfill
       \hbox{DPSU-97-11}
}}}

\vskip9mm
\begin{center}

{\large\bf 
Quasi-Hopf twistors for elliptic quantum groups} 

\vskip9mm

{\csc Michio Jimbo}$\,^{1}$,
{\csc Hitoshi Konno}$\,^{2}$,
{\csc Satoru Odake}$\,{}^3$ and
{\csc Jun'ichi Shiraishi}$\,{}^4$
\\ 
%
{\baselineskip=15pt
\it\vskip.35in 
\setcounter{footnote}{0}\renewcommand{\thefootnote}{\arabic{footnote}}
\footnote{e-mail address : jimbo@kusm.kyoto-u.ac.jp}
Division of Mathematics, Graduate School of Science,\\
Kyoto University, Kyoto 606-8502, Japan\\
\vskip.1in 
\footnote{e-mail address : konno@mis.hiroshima-u.ac.jp}
Department of Mathematics, Faculty of Integrated Arts and Sciences,\\
Hiroshima University, Higashi-Hiroshima 739-8521,  Japan\\
\vskip.1in 
\footnote{e-mail address : odake@azusa.shinshu-u.ac.jp}
Department of Physics, Faculty of Science, \\
Shinshu University, Matsumoto 390-8621, Japan\\
\vskip.1in 
\footnote{e-mail address : shiraish@momo.issp.u-tokyo.ac.jp}
Institute for Solid State Physics, \\
University of Tokyo, Tokyo 106-0032, Japan \\
}
\end{center}

\vskip6mm

\begin{center}
{\sl Dedicated to Professor Mikio Sato on the occasion of his seventieth 
birthday}
\end{center}

\begin{abstract}
The Yang-Baxter equation admits two classes of elliptic solutions, 
the vertex type and the face type.
On the basis of these solutions, two types of 
elliptic quantum groups have been introduced 
(Foda et al.\cite{FIJKMY}, Felder \cite{Fel95}). 
Fr\o nsdal \cite{Fron,Fron1} made a penetrating observation that 
both of them are quasi-Hopf algebras, 
obtained by twisting the standard quantum affine algebra $U_q(\goth{g})$. 
In this paper we present an explicit formula for the twistors 
in the form of an infinite product of the universal $R$ matrix 
of $U_q(\goth{g})$.  
We also prove the shifted cocycle condition for the twistors, 
thereby completing Fr\o nsdal's findings. 

This construction entails that, 
for generic values of the deformation parameters, 
representation theory for $U_q(\goth{g})$ carries over 
to the elliptic algebras, including such objects as 
evaluation modules, highest weight modules and vertex operators. 
In particular, we confirm the conjectures of Foda et al. 
concerning the elliptic algebra $\Aqp{}(\slth)$. 
\end{abstract}

q-alg/9712029 (to appear in Transformation Groups)

\newpage

\setcounter{section}{0}
\setcounter{equation}{0}
\section{Introduction}

\subsection{Elliptic algebras} 

Among the integrable models based on the Yang-Baxter equation (YBE), 
those related to elliptic solutions occupy a fundamental place. 
Elliptic algebras, or elliptic quantum groups, are certain 
algebraic structures introduced to account for these elliptic models. 
Nevertheless, the complexity of 
elliptic algebras has evaded their understanding for quite some time. 
The solutions of YBE ($R$-matrices) are classified into two types, 
vertex-type and face-type. 
Accordingly there are two types of elliptic algebras. 

The vertex-type elliptic algebras are associated with 
the $R$-matrix $R(u)$ of Baxter \cite{Bax72} and Belavin \cite{Bel}.
The first example of this sort 
is the Sklyanin algebra \cite{Skl82}, designed as an elliptic 
deformation of the Lie algebra $\slt$. 
(An extension to $\sln$ was discussed by Cherednik \cite{Cher85}.) 
It is presented by the `$RLL$'-relation 
\bea
R^{(12)}(u_1-u_2)L^{(1)}(u_1)L^{(2)}(u_2)
=L^{(2)}(u_2)L^{(1)}(u_1)R^{(12)}(u_1-u_2), 
\lb{vRLL}
\ena 
together with a specific choice of the form for $L(u)$. 
Here and after, the superscript $(1),(2),\cdots$ will 
refer to the tensor components. 
For further development concerning 
the Sklyanin algebra, see Feigin-Odesskii \cite{FeOd95}  
and references therein. 
An affine version of the Sklyanin algebra 
 (deformation of $\slth$) was proposed by Foda et al. \cite{FIJKMY}. 
The main point of \cite{FIJKMY} 
was to incorporate a central element $c$ 
by modifying the $RLL$-relation to 
\bea
R^{(12)}(u_1-u_2,r)L^{(1)}(u_1)L^{(2)}(u_2)
=L^{(2)}(u_2)L^{(1)}(u_1)R^{(12)}(u_1-u_2,r-c), 
\lb{vRLL2}
\ena
where $r$ denotes the elliptic modulus contained in  $R(u)=R(u,r)$. 
In both these works, the coalgebra structure was missing. 
 
The face-type algebras are based on $R$-matrices 
of Andrews, Baxter, Forrester \cite{ABF} and generalizations 
\cite{fusion,JMO1,JMO3}. 
In this case, besides the elliptic modulus, 
$R$ and $L$ depend also on extra parameter(s) $\lambda$. 
As Felder has shown \cite{Fel95}, 
the $RLL$ relation undergoes a `dynamical' shift 
by elements $h$ of the Cartan subalgebra $\goth{h}$, 
\bea
&&R^{(12)}(u_1-u_2,\la+h)L^{(1)}(u_1,\la)L^{(2)}(u_2,\la+h^{(1)})
\nonumber\\
&&\qquad=L^{(2)}(u_2,\la)L^{(1)}(u_1,\la+h^{(2)})R^{(12)}(u_1-u_2,\la). 
\lb{fRLL}
\ena
Likewise the YBE itself is modified to a dynamical one, 
see \eqref{DYBE} below. 
As we shall see, a central extension of this algebra is obtained by  
introducing further a shift of the elliptic modulus 
analogous to \eqref{vRLL2} (see \eqref{RLL1}--\eqref{RLL2}
and the remark following them). 

The face-type algebra has also been given an alternative 
formulation in terms of the Drinfeld currents. 
This is the approach adopted by Enriquez and Felder \cite{EF} 
and one of the authors \cite{Konno}.\footnote{
The central extension in \cite{Konno} is different 
from the one in \cite{EF}, and is closer to that of the present paper.}
The Drinfeld currents are suited to deal with 
infinite dimensional representations. 
We plan to discuss this subject in a separate publication. 


\subsection{Quasi-Hopf twist}

Babelon et al.\cite{BBB} have pointed out 
that the natural framework for dealing with 
dynamical YBE is Drinfeld's theory of quasi-Hopf algebras \cite{QHA}. 
Since the work \cite{BBB} is a prototype of our construction, let us recall 
their result. 
Consider the simplest quantum group $A=U_q(\slt)$ with standard 
generators $e,f,h$. 
Given an arbitrary invertible element $F\in A^{\otimes 2}$, 
we can modify the coproduct $\Delta$ and the universal $R$ matrix $\cR$ by 
\bea
&&\tilde{\Delta}(a)=F\Delta(a)F^{-1}\qquad (a\in A),
\lb{newD}\\
&&\tilde{\cR}=F^{(21)}\cR F^{-1}.
\lb{newR}
\ena
Here if $F=\sum a_i\otimes b_i$, then $F^{(21)}=\sum b_i\otimes a_i$. 
In general, the new coproduct $\tilde{\Delta}$ is no longer 
coassociative, and defines on $A$ a quasi-Hopf algebra structure. 
The new $R$ matrix $\tilde{\cR}$ satisfies a YBE-type equation, 
which is somewhat complicated (see \eqref{RPhi}). 
As Babelon et al. showed, this `twisting' procedure 
leads to an interesting result when $F=F(\la)$ depends on a parameter $\la$ 
in such a way that the shifted cocycle condition holds:
\bea
F^{(12)}(\la)(\Delta\otimes \id)F(\la)
=F^{(23)}(\la+h^{(1)})(\id \otimes \Delta)F(\la).
\lb{cocy}
\ena
If this is the case, then the YBE-type equation for 
$\tilde{\cR}=\cR(\la)$ simplifies to the dynamical YBE 
\bea
\cR^{(12)}(\la+h^{(3)})\cR^{(13)}(\la)\cR^{(23)}(\la+h^{(1)})
=
\cR^{(23)}(\la)\cR^{(13)}(\la+h^{(2)})\cR^{(12)}(\la).
\lb{DYBE}
\ena
An explicit formula for such an $F(\la)$ 
was given in \cite{BBB} as a formal power series in $q^{2\la}$. 
We shall refer to $F(\la)$ as `twistor'. 

A key observation due to Fr\o nsdal \cite{Fron1} is 
that the $RLL$ relations for the elliptic algebras 
of both types, \eqref{vRLL2} and \eqref{fRLL}, 
arise by the same mechanism as above. 
Namely, there exist two types of twistors which give rise to 
different comultiplications on the quantum affine algebras $U_q(\goth{g})$, 
and the resultant quasi-Hopf algebras 
are nothing but the two types of elliptic quantum groups. 

To substantiate this statement, we must find the corresponding 
twistors as elements in 
$ U_q(\goth{g})^{\otimes 2}$ 
satisfying the shifted cocycle condition \eqref{cocy}. 
Fr\o nsdal \cite{Fron,Fron1} launched a search 
for the twistor in the form of a formal series 
\bea
F(\la)=1+\sum_{m\ge 1}\sum_{i_1,\ldots,i_m} t_{i_1,\ldots,i_m}(\la)
e_{i_1}\ldots e_{i_m}\otimes \tau^m(f_{i_1})\ldots \tau^m(f_{i_m}), 
\lb{Fser}
\ena
where $e_i$, $f_i$ are the Chevalley generators, 
$t_{i_1,\ldots,i_m}(\la)$ are certain functions of 
the `Cartan' generators $h_i\otimes 1$ and $1\otimes h_i$,
and $\tau$ is a diagram automorphism.\footnote{
For face-type algebras $\tau=\id$, 
while for the vertex-type 
$\goth{g}=\slnh$ and $\tau$ is a cyclic diagram automorphism.
See section 2.}
Substituting \eqref{Fser} into \eqref{cocy}, 
he obtained a recursion relation that determines the 
coefficients $t_{i_1,\ldots,i_m}(\la)$ uniquely. 
Though a proof of the full cocycle condition \eqref{cocy} 
was left open, 
this construction was shown to reproduce correctly 
the classical limit \cite{Fron1} and Baxter's $R$ matrix \cite{Fron}. 
Another important observation presented 
in the work \cite{Fron1} is that the twistor 
has an infinite product form
\be
F(\la)=\cdots F^3(\la)F^2(\la)F^1(\la),
\en
and that the coefficients of each $F^m(\la)$ 
resemble those of the universal $R$ matrix of $U_q(\goth{g})$. 

We remark that 
the quasi-Hopf structure of the face-type algebra for $\slth$ 
was studied in detail by Enriquez and Felder \cite{EF} 
from a different point of view. 

\subsection{The present work}

The aim of the present article is to 
complete the works of Babelon et al. and Fr\o nsdal,  
by making explicit the aforementioned connection between the twistor and 
the universal $R$ matrix, 
and supplying a proof of the shifted cocycle condition. 
We construct the two types of 
twistors in the form of an infinite product of 
the universal $R$ matrix, 
\bea
F(\la)=\cdots \left(\varphi_\la^3\otimes\id\right)\Bigl(q^T\cR\Bigr)^{-1}
\left(\varphi_\la^2\otimes\id\right)\Bigl(q^T\cR\Bigr)^{-1}
\Bigl(\varphi_\la\otimes\id\Bigr)\Bigl(q^T\cR\Bigr)^{-1},
\lb{twistor}
\ena
where $\varphi_\la$ denotes a certain automorphism of 
$U_q(\goth{g})^{\otimes 2}$ depending on $\la$, 
and $T$ is an element of $\goth{h}\otimes\goth{h}$. 
For the face-type algebras, 
$\la$ is taken from the Cartan subalgebra, 
so that $F(\la)$ carries 
the same number of parameters as the rank of $\goth{g}$. 
When $\goth{g}$ is of affine type, the elliptic modulus appears 
as one of these parameters. 
For the vertex-type algebras, 
$\la$ is proportional to the central element.
In this case the twistor, to be denoted $E(r)$, 
depends on only one parameter $r$ which is the elliptic modulus. 

It is an old idea to construct the elliptic 
$R$ matrices and $L$ operators from the trigonometric ones 
by an `averaging' procedure over the periods 
\cite{Fa82,ReFa83,Ta85}. 
Our formula, \eqref{twistor} along with \eqref{newR}, 
may be viewed as implementing this idea 
at the level of the universal $R$ matrix. 

Following \cite{FIJKMY}, we shall denote the 
quasi-Hopf algebras associated with the 
vertex-type twistor $E(r)$ 
by the symbol $\Aqp{}(\goth{g})$ (where $\goth{g}=\slnh$), 
and the one associated with the face-type twistor $F(\la)$ 
by $\Bqla(\goth{g})$. 
As algebras they are the same as 
the underlying $U_q(\goth{g})$. 
Hence the representation theory for them should stay the same. 
Strictly speaking, the twistors are 
only formal power series with coefficients in $U_q(\goth{g})$, 
but we expect they make sense in `good' category of representations
and for generic values of the parameters. 
In such a case, the whole representation theory
including evaluation and highest weight modules and 
vertex operators carry over to the elliptic algebras. 
We derive the commutation relations of vertex operators and the 
intertwining relations, regarding them as formal power series. 
In the special case of $\Aqp{}(\slth)$ we recover the 
formulas conjectured in \cite{FIJKMY}. 
As Fr\o nsdal has shown \cite{Fron,Fron2} for $\Aqp{}(\slth)$, 
the formal series for the twistors in  evaluation modules
converge and can be computed explicitly.
On the other hand, it is a non-trivial problem to 
compute their images in highest weight representations.
We hope to come back to this issue in the future. 

The text is organized as follows. 

In section 2, after preparing the notation, 
we present the formulas for the twistors. 
We then give a proof of the shifted cocycle condition. 
In Section 3, we discuss the examples 
$\Bqla(\slt)$, $\Bqla(\slth)$ and $\Aqp{}(\slth)$, 
and compute the images of the twistor and the 
universal $R$ matrix in the two-dimensional 
evaluation representation. 
In section 4, we define $L$ operators and vertex operators 
for the elliptic algebras out of those of $U_q(\goth{g})$, 
and derive various commutation relations among them. 
In particular we derive a relation between the $L^+$ and $L^-$ operators 
proposed earlier in \cite{FIJKMY}
(see \eqref{LpLm}, \eqref{Lrel}). 
In Appendix we review the basics of quasi-Hopf algebras.


\setcounter{section}{1}
\setcounter{equation}{0}
\section{Quasi-Hopf twistors}

In this section we construct the twistors 
which give rise to the elliptic algebras. 
For the face type algebras, the twistor 
$F=F(\la)\in U^{\otimes2}$ depends on a parameter $\la$
running over the Cartan subalgebra $\goth{h}$. 
For the vertex type algebras, the twistor 
$E(r)$ depends on a single parameter $r\in\C$. 
Both of them are solutions of the shifted cocycle condition 
(see \eqref{facecocy}, \eqref{vertexcocy} below). 

\subsection{Quantum groups}\lb{sec:2.1}

First let us fix the notation.
Let $\goth{g}$ be the Kac-Moody Lie algebra associated with 
a symmetrizable generalized Cartan matrix $A=({a_{ij}})_{i,j\in I}$. 
We fix an invariant inner product $(~,~)$ on 
the Cartan subalgebra $\goth{h}$
and identify $\goth{h}^*$ with $\goth{h}$ via $(~,~)$.  
If $\{\alpha_i\}_{i\in I}$ denotes the set of 
simple roots, then $(\alpha_i,\alpha_j)=d_ia_{ij}$, where 
$d_i={1\over 2}(\alpha_i,\alpha_i)$. 

Consider the corresponding quantum group $U=U_q(\goth{g})$. 
For simplicity of presentation, 
we choose to work over the ground ring $\C[[\hbar]]$ with $q=e^\hbar$. 
The algebra $U$ has generators $e_i$, $f_i$ ($i\in I$) 
and $h$ ($h\in \goth{h}$), 
satisfying the standard relations 
\bea
&&[h,h']=0\qquad (h,h'\in \goth{h}),
\lb{Uqrel0}\\
&&
[h,e_i]=(h,\alpha_i) e_i, \qquad [h,f_i]=-(h,\alpha_i) f_i
\qquad (i\in I, h\in \goth{h}),
\lb{Uqrel1}\\
&&[e_i,f_j]=\delta_{ij}\frac{t_i-t_i^{-1}}{q_i-q^{-1}_i}
\qquad (i,j\in I),
\lb{Uqrel2}
\ena
and the Serre relations which we omit. 
In \eqref{Uqrel2} we have set $q_i=q^{d_i}$, $t_i=q^{\alpha_i}$. 
We adopt the Hopf algebra structure given as follows. 
\bea
&&\Delta(h)=h\otimes 1+1 \otimes h,\\
&&
\Delta(e_i)=e_i\otimes 1 + t_i\otimes e_i,
\quad
\Delta(f_i)=f_i\otimes t_i^{-1}+1\otimes f_i, 
\label{copro}\\
&&\ve(e_i)=\ve(f_i)=\ve(h)=0, 
\label{counit}\\
&&S(e_i)=- t_i^{-1} e_i,\quad 
  S(f_i)=- f_i t_i,\quad 
  S(h)=- h,
\label{antipode} 
\ena
where $i\in I$ and $h\in \goth{h}$. 

Let $\cR\in U^{\otimes 2}$ denote the universal $R$ matrix of $U$. 
It has the form 
\begin{eqnarray} 
&&\cR=q^{-T} \cC ,
\\
&&\cC=\sum_{\beta\in Q^+} q^{(\beta,\beta)}
\Bigl(q^{-\beta}\otimes q^{\beta}\Bigr)\cC_\beta
\nn\\
&&\phantom{\cC}=
1-\sum_{i\in I}(q_i-q_i^{-1})e_i t_i^{-1}\otimes t_i f_i+\cdots.
\end{eqnarray}
Here the notation is as follows. 
Take  a basis $\{h_l\}$ of $\goth{h}$, and its dual basis $\{h^l\}$. 
Then 
\bea
T=\sum_l h_l\otimes h^l 
\lb{T}
\ena
denotes the canonical element of $\goth{h}\otimes \goth{h}$. 
The element 
$\cC_\beta=\sum_j u_{\beta,j}\otimes u^j_{-\beta}$ 
is the canonical element of $U^+_\beta\otimes U^-_{-\beta}$ 
with respect to a certain Hopf pairing, 
where $U^+$ (resp. $U^-$) denotes the subalgebra of $U$ 
generated by the $e_i$ (resp. $f_i$), and 
$U^\pm_{\pm \beta}$ ($\beta\in Q^+$)
signifies the homogeneous components 
with respect to the natural gradation by $Q^+=\sum_i\Z_{\ge 0}\alpha_i$.
(For the details the reader is referred e.g. to \cite{Dri86,Tani}.)
We shall need the following 
basic properties of the universal $R$ matrix:
\begin{eqnarray}
&&\Delta'(a)=\cR \Delta(a)\cR^{-1} \qquad \forall a\in U, 
\label{qt1}\\
&&
\left(\Delta\otimes \id\right)\cR
=\cR^{(13)}\cR^{(23)},
\label{qt2}\\
&&
\left(\id\otimes \Delta\right)\cR
=\cR^{(13)}\cR^{(12)},
\label{qt3}\\
&&\left(\varepsilon\otimes\id\right)\cR
=\left(\id\otimes\varepsilon\right)\cR
=1.
\label{qt4}
\end{eqnarray}
Here $\Delta'=\sigma\circ \Delta $ signifies the opposite coproduct, 
$\sigma$ being the flip of the 
tensor components $\sigma(a\otimes b)=b\otimes a$.
{}From (\ref{qt1})--(\ref{qt3}) follows the Yang-Baxter equation 
\begin{eqnarray}
&&\cR^{(12)} \cR^{(13)}\cR^{(23)}
=\cR^{(23)}\cR^{(13)}\cR^{(12)}. 
\label{yb}
\end{eqnarray}

\subsection{Face type twistors}\lb{sec:2.2}

We are now in a position to describe the twistors for 
face type elliptic algebras. 
Let us prepare some notation. 

Let $\rho\in\goth{h}$ be an element such that 
$(\rho,\alpha_i)=d_i$ for all $i\in I$. 
Let $\phi$ be an automorphism of $U$ given by 
\bea
&&\phi=\Ad(q^{\frac{1}{2}\sum_l h_lh^l-\rho}), 
\lb{phi}
\ena
where $\{h_l\}$, $\{h^l\}$ are as in \eqref{T}. 
In other words, 
\be
&&\phi(e_i)=e_i t_i,\qquad \phi(f_i)=t_i^{- 1}f_i,\qquad
\phi(q^{h})=q^{h}.
\en
Since 
\bea
\Ad(q^T)\circ(\phi\otimes\phi)=
\Ad(q^{\frac{1}{2}\sum_l \Delta(h_lh^l)-\Delta(\rho)}),
\lb{phiphi}
\ena
we have 
\bea
&&\Ad(q^{T})\circ(\phi\otimes\phi)\circ \Delta
=\Delta \circ \phi. 
\lb{delphi}
\ena
For $\la\in \goth{h}$, introduce an automorphism 
\bea
\varphi_\la=\Ad(q^{\sum_l h_lh^l+2(\la-\rho)})
=\phi^2\circ\Ad(q^{2\la}).
\lb{phila}
\ena
Then the expression 
\bea
\Bigl(\varphi_\la\otimes\id\Bigr)\Bigl(q^T\cR\Bigr)
\lb{qTR1}
\ena
is a formal power series in the variables 
$x_i=q^{2(\la,\alpha_i)}$ ($i\in I$) of the form 
\be
1-\sum_i(q_i-q_i^{-1})x_i e_it_i\otimes t_i f_i+\cdots. 
\en
We define the twistor $F(\la)$ as follows.
\begin{dfn}[Face type twistor]\lb{dfn:face}
\bea
F(\la)&=&\cdots 
\Bigl(\varphi_\la^2\otimes\id\Bigr)\Bigl(q^T\cR\Bigr)^{-1}
\Bigl(\varphi_\la\otimes\id\Bigr)\Bigl(q^T\cR\Bigr)^{-1}
\nn\\
\phantom{F(\la)}&=&
\mathop{\prod_{k\geq 1}}^{\kuru}
\Bigl(\varphi_\la^k\otimes\id\Bigr)\Bigl(q^T\cR\Bigr)^{-1}.
\label{facetwistor}
\ena
\end{dfn}
Here and after, we use the ordered product symbol 
$\displaystyle\mathop{\prod_{k\geq 1}}^{\kuru}A_k=\cdots A_3A_2A_1$. 
Note that the $k$-th factor in the product \eqref{facetwistor}
is a formal power series 
in the $x_i^k$ with leading term $1$, 
and hence the infinite product makes sense. 
We shall refer to \eqref{facetwistor} as a {\it face type} twistor. 

Our main result is the following. 
\begin{thm}\lb{thm:face}
The twistor \eqref{facetwistor} satisfies the
shifted cocycle condition 
\begin{eqnarray}
F^{(12)}(\la)(\Delta\otimes \id)F(\la)
=F^{(23)}(\la+h^{(1)})(\id \otimes \Delta)F(\la).
\lb{facecocy}
\end{eqnarray}
We have in addition
\bea
&&\left(\varepsilon\otimes\id\right)F(\la)
=\left(\id\otimes\varepsilon\right)F(\la)=1.
\lb{epF}
\ena
\end{thm}
A proof of Theorem \ref{thm:face}
will be given in subsection \ref{sec:2.4}. 
In \eqref{facecocy}, if $\la=\sum_l\la_lh^l$, 
then $\la+h^{(1)}$ means $\sum_l(\la_l+h_l^{(1)})h^l$. 
Hence we have, for example, 
\be
&&\Ad(q^{2l T^{(12)}}) F^{(23)}_k(\la)=F^{(23)}_k(\la+ l h^{(1)}),
\\
&&\Ad(q^{2l T^{(13)}}) F^{(23)}_k(\la)=F^{(23)}_k(\la- l h^{(1)}).
\en

For convenience, let us give a name to the quasi-Hopf algebra
associated with the twistor \eqref{facetwistor}. 
As for the generalities on quasi-Hopf algebras, see Appendix \ref{app:a}.

\begin{dfn}[Face type algebra]\lb{dfn:Bqla} 
We define the quasi-Hopf algebra $\Bqla(\goth{g})$ of face type
to be the set 
$(U_q(\goth{g}),\Delta_\la,\ve,\Phi(\la),\cR(\la))$
together with 
$\alpha_\la=\sum_i S(d_i)e_i$, 
$\beta_\la=\sum_i f_i S(g_i)$ and 
the antiautomorphism $S$ defined by (\ref{antipode}).
Here $\ve$ is defined by (\ref{counit}),
\bea
&&\Delta_{\la}(a)=F^{(12)}(\la) \,\Delta(a)\,  F^{(12)}(\la)^{-1},
\lb{facecopro}\\
&&\cR(\la)=F^{(21)}(\la)\,\cR\, F^{(12)}(\la)^{-1},
\lb{faceR}\\
&&\Phi(\la)=F^{(23)}(\la)F^{(23)}(\la+h^{(1)})^{-1},
\lb{facephi}
\ena
and $\sum_i d_i\otimes e_i=F(\la)^{-1}$, $\sum_i f_i \otimes g_i=F(\la)$.
\end{dfn}

Let us consider the case where $\goth{g}$ is of affine type, 
in which we are mainly interested. 
Let $c$ be the canonical central element and 
$d$ the scaling element. 
We set 
\be
\la-\rho=rd+s'c+\bar{\la}-\bar{\rho}
\qquad (r,s'\in \C),
\en
where $\bar{\la}$ stands for the classical part of $\la\in\goth{h}$. 
Denote by $\{\bar{h}_j\}$, $\{\bar{h}^j\}$ the 
classical part of the dual basis of $\goth{h}$. 
Since $c$ is central, $\varphi_\la$ is independent of $s'$. 
Writing $p=q^{2r}$, 
we have 
\be
\varphi_\la=\Ad(p^dq^{2cd})\circ\bar{\varphi}_\la,
\qquad 
\bar{\varphi}_\la=\Ad(q^{\sum \bar{h}_j\bar{h}^j+2(\bar{\la}-\bar{\rho})}).
\en
Set further\footnote{The notation $\cR(z)$ conflicts that of 
$\cR(\la)$. Hopefully there is no confusion.} 
\bea
&&\cR(z)=\Ad(z^d\otimes 1)(\cR),
\lb{Rg}
\\
&&F(z,\la)=\Ad(z^d\otimes 1)(F(\la)),
\lb{Fg}
\\
&&\cR(z,\la)=\Ad(z^d\otimes 1)(\cR(\la))
=\sigma(F(z^{-1},\la))\cR(z)F(z,\la)^{-1}.
\lb{Rg2}
\ena
Here $\sigma$ denotes the flip of the tensor components. 
\eqref{Rg}, \eqref{Fg} are formal power series in $z$, whereas 
\eqref{Rg2} contains both positive and negative powers of $z$. 
Note that $q^{c\otimes d+d\otimes c}\cR(z)\bigl|_{z=0}$ reduces to the 
universal $R$ matrix of $U_q(\bar{\goth{g}})$ 
corresponding to the underlying finite dimensional Lie algebra
$\bar\goth{g}$. 
{}From the definition \eqref{facetwistor} of $F(\la)$ we have the 
difference equation 
\bea
&&F(pq^{2c^{(1)}}z,\la)=(\bar{\varphi}_\la\otimes\id)^{-1}
\bigl(F(z,\la)\bigr)\cdot q^{T}\cR(pq^{2c^{(1)}}z),
\lb{Fdiff}
\\
&&F(0,\la)=F_{\bar{\goth{g}}}(\bar{\la}),
\lb{Fini}
\ena
where $F_{\bar{\goth{g}}}(\bar{\la})$ signifies the twistor 
corresponding to $\bar\goth{g}$. 
\medskip

\subsection{Vertex type twistors}\lb{sec:2.3} 

When $\goth{g}=\slnh$, 
it is possible to construct a different type of twistor. 
We call it {\it vertex type}. 
In this subsection, $U$ will denote $U_q(\slnh)$. 

Let us write $h_i=\alpha_i$ ($i=0,\ldots,n-1$). 
A basis of $\goth{h}$ is 
$\{h_0,\ldots,h_{n-1},d\}$. 
The element $d$ gives the homogeneous grading, 
\be
[d,e_i]=\delta_{i0}e_i,
\quad
[d,f_i]=-\delta_{i0}f_i,
\en
for all $i=0,\ldots,n-1$. 
Let the dual basis be $\{\La_0,\ldots,\La_{n-1},c\}$. 
The $\La_i$ are the fundamental weights 
and $c$ is the canonical central element.
Let $\tau$ be the automorphism of $U$ such that 
\be
\tau(e_i)=e_{i+1\smod n},
\qquad
\tau(f_i)=f_{i+1\smod n},
\qquad
\tau(h_i)=h_{i+1\smod n}
\en
and $\tau^n=\id$. 
Then we have 
\be
\tau(\La_i)=\La_{i+1\smod n}-\frac{n-1-2i}{2n}c.
\en
The element 
$\rho=\sum_{i=0}^{n-1}\La_i$ is invariant under $\tau$. 
It gives the principal grading 
\be
[\rho,e_i]=e_i,
\quad
[\rho,f_i]=-f_i,
\en
for all $i=0,\ldots,n-1$. 
Note also that 
\be
&&(\tau\otimes\tau)\circ \Delta=\Delta\circ\tau,
\\
&&(\tau\otimes\tau) (\cC_\beta)=\cC_{\tau(\beta)}.
\en

For $r\in \C$, we introduce an automorphism 
\bea
\widetilde{\varphi}_r=\tau\circ\Ad\Bigl(q^{\frac{2(r+c)}{n}\rho}\Bigr).
\lb{philatilde}
\ena
Here and after, quantities related to the vertex type algebras will be 
denoted with the symbol $\widetilde{\phantom{\varphi}}$. 
Set 
\be 
\widetilde{T}=\frac{1}{n}\Bigl(\rho\otimes c+c\otimes \rho
-\frac{n^2-1}{12}c\otimes c\Bigr).
\en
Then 
\bea
\Bigl(\widetilde{\varphi}_r\otimes\id\Bigr)
\Bigl(q^{\widetilde{T}}\cR\Bigr)^{-1}
\lb{qTR2}
\ena
is a formal power series in $p^{\frac{1}{n}}$ where $p=q^{2r}$. 
Unlike the previous case of \eqref{qTR1}, \eqref{qTR2} 
is a formal series with a non-trivial leading term 
$q^{T-\widetilde{T}}\left(1+\cdots\right)$. 
Nevertheless, the $n$-fold product 
\be
\mathop{\prod_{n\ge k\ge 1}}^{\kuru}
\Bigl(\widetilde{\varphi}_r^k\otimes\id\Bigr)
\Bigl(q^{\widetilde{T}}\cR\Bigr)^{-1}
\en
takes the form $1+\cdots$, because of the relation
\be
\sum_{k=1}^n\left(\tau^k\otimes\id\right)\left(T-\widetilde{T}\right)=0.
\en
We now define the vertex type twistor $E(r)$ as follows.
\begin{dfn}[Vertex type twistor]\lb{dfn:vertex}
\bea
E(r)&=&
\mathop{\prod_{k\geq 1}}^{\kuru}
\left(\widetilde{\varphi}_r^k\otimes\id\right)
\Bigl(q^{\widetilde{T}}\cR\Bigr)^{-1}.
\label{vertextwistor}
\ena
\end{dfn}
The infinite product $\displaystyle\mathop{\prod_{k\geq 1}}^{\kuru}{}$
is to be understood as
$\displaystyle\lim_{N\rightarrow\infty}
\mathop{\prod_{nN\ge k\ge 1}}^{\kuru}{}$. 
In view of the remark made above, $E(r)$ is 
a well defined formal series in $p^{\frac{1}{n}}$. 

\begin{thm}\lb{thm:vertex}
The twistor \eqref{vertextwistor} satisfies the 
shifted cocycle condition 
\begin{eqnarray}
E^{(12)}(r)(\Delta\otimes \id)E(r)
=E^{(23)}(r+c^{(1)})(\id \otimes \Delta)E(r).
\lb{vertexcocy}
\end{eqnarray}
We have in addition
\bea
&&\left(\varepsilon\otimes\id\right)E(r)
=\left(\id\otimes\varepsilon\right)E(r)=1.
\lb{epE}
\ena
\end{thm}

\begin{dfn}[Vertex type algebra]\lb{dfn:Aqp} 
We define the quasi-Hopf algebra $\Aqp{}(\slnh)$ 
($p=q^{2r}$) of vertex type to be the set 
$(U_q(\goth{g}),\Delta_r,\ve,\Phi(r),\cR(r))$ 
together with $\alpha_r=\sum_i S(d_i)e_i$, 
$\beta_r=\sum_i f_i S(g_i)$ and 
the antiautomorphism $S$ defined by (\ref{antipode}). 
Here $\ve$ is defined by (\ref{counit}),
\bea
&&\Delta_{r}(a)=E^{(12)}(r) \,\Delta(a)\,  E^{(12)}(r)^{-1},
\lb{vertexcopro}\\
&&\cR(r)=E^{(21)}(r)\,\cR\, E^{(12)}(r)^{-1},
\lb{vertexR}\\
&&\Phi(r)=E^{(23)}(r)E^{(23)}(r+c^{(1)})^{-1},
\lb{vertexphi}
\ena
and
$\sum_i d_i\otimes e_i=E(r)^{-1}$,
$\sum_i f_i \otimes g_i=E(r)$.
\end{dfn}

Let us set 
\be
&&\widetilde{\cR}'(\ze)=\left(\Ad(\ze^\rho)\otimes\id\right)
(q^{\widetilde{T}}\cR),
\\
&&
E(\ze,r)=\left(\Ad(\ze^\rho)\otimes\id\right)E(r).
\en
In just the same way as in the face type case, 
the definition \eqref{vertextwistor} can be alternatively described as the 
unique solution of the difference equation 
\be
&&E(p^{1/n}q^{2c^{(1)}/n}\ze,r)=
(\tau\otimes\id)^{-1}\bigl(E(\ze,r)\bigr)\cdot 
\widetilde{\cR}'(p^{1/n}q^{2c^{(1)}/n}\ze),
\lb{Ediff}
\\
&&E(0,r)=1,
\lb{Eini}
\en
where $p=q^{2r}$.

\subsection{Proof of the shifted cocycle condition}\lb{sec:2.4}

Let us prove the shifted cocycle condition \eqref{facecocy}
for the face type twistors. 
For $k=0,1,\cdots$ we set 
\bea
&&F_k(\la)=(\phi^{2k}\otimes \id) \cC(\la)^{-1},
\lb{Fk}\\
&&
\cC(\la)=\Ad(q^{2\la} \otimes 1) (q^T \cR)
=\Ad(1\otimes q^{-2\la}) (q^T \cR).
\lb{Cla}
\ena
Then the twistor \eqref{facetwistor} can be written as 
\be
F(\la)=\mathop{\prod_{k\geq 1}}^{\kuru}F_k(k\la).
\en
We have the invariance $[\Delta(h),F_k(\la)]=0$, and in particular 
\bea
&&\Ad(q^T)\circ(\phi\otimes\phi)\left(F_k(\la)\right)=F_k(\la). 
\lb{Finv}
\ena

{}From the properties (\ref{qt1})--(\ref{qt3})
of the universal $R$ matrix, we find
\begin{lem}
\bea
&&(\Delta \otimes \id) \cC(\la)= 
 \cC^{(13)}(\la-\frac{1}{2}h^{(2)})  \cC^{(23)}(\la),
\label{c1}\\
&&(\id\otimes \Delta) \cC(\la)= 
\cC^{(13)}(\la+\frac{1}{2}h^{(2)})  \cC^{(12)}(\la),
\label{c2}\\
&&    \cC^{(12)}(\la)
   \cC^{(13)}(\la+\mu-\frac{1}{2}h^{(2)}) 
   \cC^{(23)}(\mu)
\nn \\
&&\qquad =
   \cC^{(23)}(\mu)
   \cC^{(13)}(\la+\mu+\frac{1}{2}h^{(2)})\cC^{(12)}(\la).
\label{c3}
\ena
\end{lem}

\begin{lem}
\begin{eqnarray}
&&(\Delta \otimes \id) F_k(\la)= 
  F_k^{(23)}(\la+kh^{(1)})F_k^{(13)}(\la+(k-\frac{1}{2})h^{(2)}),
\label{f1}\\
&&(\id\otimes \Delta) F_k(\la)= 
  F_k^{(12)}(\la)F_k^{(13)}(\la+\frac{1}{2}h^{(2)}),
\label{f2}\\
&& F_k^{(12)}(\la) 
   F_{k+l}^{(13)}(\la+\mu+(l+\frac{1}{2})h^{(2)}) 
   F_l^{(23)}(\mu+lh^{(1)}) \nn \\
&&\qquad =
   F_l^{(23)}(\mu+lh^{(1)}) 
   F_{k+l}^{(13)}(\la+\mu+(l-\frac{1}{2})h^{(2)}) 
   F_k^{(12)}(\la).
\label{f3}
\end{eqnarray}
\end{lem}
\proof
Using (\ref{c1}) we have 
\be
\mbox{LHS of }(\ref{f1})\!\!\!\!\!\!\!\!
&&=(\Delta\otimes \id) (\phi^{2k}\otimes \id) \cC(\la)^{-1} \\
&&= \Ad(q^{2k T^{(12)}})
(\phi^{2k}\otimes\phi^{2k}\otimes \id)
\left( \cC^{(23)}(\la)^{-1} \cC^{(13)}(\la-\frac{1}{2}h^{(2)})^{-1}\right) 
\\
&&=F_k^{(23)}(\la+kh^{(1)}) F_k^{(13)}(\la+(k-\frac{1}{2})h^{(2)}). 
\en
In the second line we used \eqref{delphi}. 
Eq.\eqref{f2} can be verified in a similar way. 
Finally, \eqref{f3} follows by applying 
$(\phi^{2(k+l)}\otimes \phi^{2l}\otimes \id)\Ad(q^{2lT^{(12)}})$ 
to \eqref{c3} and noting \eqref{Finv}. 
\qed

\begin{lem}
For $l\in \Z_{\geq 0}$, we have the equality 
\begin{eqnarray}
&&\mathop{\prod_{l\ge k\ge 1}}^{\kuru}F_k^{(23)}(k\la+kh^{(1)})
\cdot (\id\otimes \Delta ) F(\la) \nn\\
&&\quad =
\mathop{\prod_{k\ge 1}}^{\kuru}
F^{(12)}_k(k\la)F^{(13)}_{k+l}((k+l)\la+(l+\frac{1}{2})h^{(2)}) \nn\\
&&\qquad \times
\mathop{\prod_{l\ge k\ge 1}}^{\kuru}
F^{(23)}_k(k\la+k h^{(1)})
F^{(13)}_k(k\la+(k-\frac{1}{2})h^{(2)}).
\label{lemma3}
\end{eqnarray}
\end{lem}
\proof
We prove \eqref{lemma3} by induction on $l$. 
The statement holds for $l=0$, 
since we have from (\ref{f2})
\begin{eqnarray*}
&&(\id\otimes \Delta)F(\la)=
\mathop{\prod_{k\ge 1}}^{\kuru}
F^{(12)}_k(k\la)
F^{(13)}_k(k\la+\frac{1}{2}h^{(2)}).
\end{eqnarray*}

Suppose the statement is correct for $l-1$. 
Then from (\ref{f3}) we obtain 
\be
\!\!\!\!\!\!
&&\phantom{=}F_l^{(23)}(l\la+l h^{(1)}) \cdot
\mathop{\prod_{k\ge 1}}^{\kuru}
F^{(12)}_k(k\la)F^{(13)}_{k+l-1}((k+l-1)\la+(l-\frac{1}{2})h^{(2)}) 
\\
\!\!\!\!\!\!&&
=F_l^{(23)}(l\la+l h^{(1)}) \cdot
\mathop{\prod_{k\ge 1}}^{\kuru}
F^{(13)}_{k+l}((k+l)\la+(l-\frac{1}{2})h^{(2)})
F^{(12)}_k(k\la)
\times F_l^{(13)}(l\la+(l-\frac{1}{2}) h^{(2)})\\
\!\!\!\!\!\!&&=
\mathop{\prod_{k\ge 1}}^{\kuru}
F^{(12)}_k(k\la)
F^{(13)}_{k+l}((k+l)\la+(l+\frac{1}{2})h^{(2)})
\times F_l^{(23)}(l\la+l h^{(1)}) 
F_l^{(13)}(l\la+(l-\frac{1}{2}) h^{(2)}). 
\end{eqnarray*}
This means that the statement holds also for $l$. 
\qed

\noindent
{\it Proof of Theorem \ref{thm:face}.}
Let $l\rightarrow\infty$ in (\ref{lemma3}). 
Then 
\be
&&\phantom{=}F^{(23)}(\la+h^{(1)}) (\id\otimes \Delta)F(\la)\\
&&=
\mathop{\prod_{k\ge 1}}^{\kuru}
F^{(12)}_k(k\la) \cdot
\mathop{\prod_{k \geq 1}}^{\kuru}
F^{(23)}_k(k\la+k h^{(1)})
F^{(13)}_k(k\la+(k-\frac{1}{2})h^{(2)})\\
&&=
F^{(12)}(\la) (\Delta\otimes \id)F(\la).
\en
The last step is from (\ref{f1}). 
The statement \eqref{epF} is evident from \eqref{qt4}. 
\qed

The case of vertex type twistors \eqref{vertexcocy} can be treated 
in an analogous manner. 
In place of the automorphism \eqref{phi}, we set 
\be
\widetilde{\phi}=\tau\circ \Ad(q^{\frac{2c}{n}\rho}).
\en
Then we have
\be
&&\Ad\Bigl(q^{2\widetilde{T}}\Bigr)\circ
(\widetilde{\phi}\otimes\widetilde{\phi})\circ\Delta
=\Delta\circ\widetilde{\phi}.
\en
The twistor can be written as 
\be
E(r)=\mathop{\prod_{k\geq 1}}^{\kuru}E_k(kr), 
\en
with the definition 
\bea
&&E_k(r)=(\widetilde{\phi}^{k}\otimes \id) \widetilde{\cC}(r)^{-1},
\lb{Ek}\\
&&
\widetilde{\cC}(r)=\Ad(q^{\frac{2r}{n}\rho} \otimes 1)
(q^{\widetilde{T}}\cR)
=\Ad(1\otimes q^{-\frac{2r}{n}\rho}) (q^{\widetilde{T}}\cR).
\lb{Cr}
\ena
We have also 
\be
&&\Ad\Bigl(q^{2\widetilde{T}}\Bigr)\circ
(\widetilde{\phi}\otimes\widetilde{\phi})
\left(E_k(r)\right)
=
E_k(r).
\en
The rest of the proof is much the same with that of 
the face type, so we omit the details. 


\setcounter{section}{2}
\setcounter{equation}{0}
\section{Examples} 

\subsection{The case $\Bqla(\sltbig)$}

Let $\goth{g}=\slt$, with the generators $e,f,h$ as in 
\eqref{Uqrel1},\eqref{Uqrel2}. 
In this case the universal $R$ matrix is given by \cite{Dri86}
\bea
&&\cR=q^{-T}\exp_{q^{2}}\left(-(q-q^{-1}) et^{-1}\otimes tf\right),
\qquad T=\frac{1}{2}h\otimes h. 
\lb{Rslt}
\ena
Here the $q$-exponential symbol is defined by 
\be
&&\exp_q(x)=\sum_{n=0}^\infty {x^n \over (n)_q!},
\qquad 
\exp_q(x)\exp_{q^{-1}}(-x)=1,
\\
&&(n)_q!=\frac{(q;q)_n}{(1-q)^n}, 
\qquad (a;q)_n=\prod_{k=0}^{n-1}(1-a q^k). 
\en
Let us set 
\be
\la=(s+1)\frac{1}{2}h, \quad w=q^{2s}, 
\en
and write $F(w)$ for $F(\la)$. 
Since 
\be
\varphi_\la=\Ad(q^{\frac{1}{2}h^2}w^{\frac{1}{2}h})
\en
and $\varphi_\la^k(e)=(q^2w)^ket^{2k}$, the formula for the twistor 
\eqref{facetwistor} becomes 
\bea
&&F(w)=\mathop{\prod_{k\geq 1}}^{\kuru}
\exp_{q^{-2}}\left((q-q^{-1})(q^2w)^k \cdot e t^{2k-1} \otimes tf\right).
\lb{Fslt1} 
\ena
Using the formula 
\be
\sum_{n=0}^\infty\frac{1}{(n)_q!(a;q)_n}b^n
=\mathop{\prod_{k\geq 0}}^{\kuru}\exp_q(a^k b)
\qquad
\mbox{ if~~ $ba=q^2ab$}, 
\en
we find 
\bea
F(w)=\sum_{n=0}^{\infty}
\frac{(q^2w)^n(q-q^{-1})^n}{(n)_{q^{-2}}!(q^{-2}w(t^2\otimes 1);q^{-2})_n}
(et)^n\otimes (tf)^n. 
\lb{Fslt2}
\ena
The formulas \eqref{Fslt1}--\eqref{Fslt2} are due to \cite{BBB,Fron,Fron1}.
In the two-dimensional representation $(\pi,\C^2)$ 
\bea
\pi(e)=E_{12}, \quad
\pi(f)=E_{21}, \quad
\pi(h)=E_{11}-E_{22},
\lb{2jigen}
\ena
with $E_{ij}$ denoting the matrix with $1$ at the $(i,j)$-th place 
and $0$ elsewhere,
we have 
\bea
(\pi\otimes\pi)F(w)=1+(q-q^{-1})\frac{w}{1-w}E_{12}\otimes E_{21},
\lb{Fslt3}
\ena

\subsection{The case $\Bqla(\slthbig)$} 

Next consider the case of the affine Lie algebra $\goth{g}=\slth$. 
Taking a basis $\{c,d,h_1\}$ of $\goth{h}$, we write 
\be
\la=(r+2)d+s'c+(s+1)\frac{1}{2}h_1 
\qquad (r,s',s\in \C).
\en
$\varphi_\la$ is independent of $s'$. 
Writing 
\bea
p=q^{2r},\quad w=q^{2s}, 
\lb{pw}
\ena
we set 
\bea
&&\cR(z)=\Ad(z^d\otimes 1)(\cR),
\lb{Rslth1}
\\
&&F(z;p,w)=\Ad(z^d\otimes 1)(F(\la)),
\lb{Fslth}
\\
&&\cR(z;p,w)=\Ad(z^d\otimes 1)(\cR(\la))
=\sigma(F(z^{-1};p,w))\cR(z)F(z;p,w)^{-1}.
\lb{Rslth2}
\ena
In particular, for $z=0$, $q^{c\otimes d+d\otimes c}\cR(0)$ 
reduces to the universal $R$ matrix \eqref{Rslt} of $U_q(\slt)$. 
{}From \eqref{Fdiff},\eqref{Fini} we have 
\bea
&&
F(pq^{2c^{(1)}}z;p,w)=(\bar{\varphi}_w^{-1}\otimes\id)
\left(F(z;p,w)\right)\times q^T\cR(pq^{2c^{(1)}}z),
\lb{facediff}
\\
&&F(0;p,w)=F_{\goth{\slt}}(w),
\lb{faceini}
\ena 
where $\bar{\varphi}_w=\Ad(q^{h_1^2/2}w^{h_1/2})$. 
In \eqref{faceini}, the right hand side means 
the twistor \eqref{Fslt1} in the previous example.

Let us calculate the image of \eqref{Fslth},\eqref{Rslth2} 
in the two-dimensional representation $(\pi,V)$, $V=\C^2$, 
where $e_1,f_1,h_1$ are mapped as in \eqref{2jigen} and 
$\pi(e_0)=\pi(f_1)$, $\pi(f_0)=\pi(e_1)$, $\pi(h_0)=-\pi(h_1)$.
We set 
\be
&&F_{VV}(z;p,w)=(\pi\otimes\pi)F(z;p,w),
\\
&&
R_{VV}(z;p,w)=(\pi\otimes\pi)\cR(z;p,w).
\en
The image 
$R_{VV}(z)=(\pi\otimes\pi)\cR(z)$ 
is known to be given as follows (see e.g.\cite{IIJMNT}). 
\bea
&&
R_{VV}(z)=\rho(z)\bR_{VV}(z),
\nn\\
&&
\rho(z)=q^{-\frac{1}{2}}\frac{(q^2z;q^4)_\infty^2}
{(z;q^4)_\infty(q^4z;q^4)_\infty},
\lb{brho}\\
&&\bR_{VV}(z)=\left(
\begin{array}{cccc}
1 & & & \\
 & b(z) & c(z) & \\
 & zc(z) & b(z) & \\
 &  &  & 1 \\
\end{array}\right),
\lb{bR1}\\
&&
b(z)=\frac{(1-z)q}{1-q^2z},
\qquad 
c(z)=\frac{1-q^2}{1-q^2z}.
\lb{bR2}
\ena
Here and after, we use the infinite product symbol
\be
&&
(z;t_1,t_2,\cdots,t_n)_\infty =\prod_{i_1,i_2,\cdots,i_n=0}^\infty
(1-z t_1^{i_1} t_2^{i_2} \cdots t_n^{i_n} ). 
\en

{}Eq. \eqref{facediff}
implies a difference equation for $F_{VV}(z;p,w)$. 
Noting $\pi(c)=0$ and 
$\pi\circ\bar{\varphi}_w=\Ad(D_w)^{-1}\circ\pi$ 
where $D_w=\mbox{diag}(1,w)$, we find 
\be
F_{VV}(pz;p,w)^{t}=R_{VV}(pz)^{t}\,K(D_w\otimes 1)^{-1}
\cdot F_{VV}(z;p,w)^{t}\,(D_w\otimes 1),
\en
where $X^t$ means the transpose of $X$, and  
we have set $K=\mbox{diag}(q^{1/2},q^{-1/2},q^{-1/2},q^{1/2})$. 
This means that each column of $F_{VV}(z;p,w)^{t}$ satisfies 
a difference equation of the same sort as the $q$-KZ equation. 
Solving this with the initial condition
which follows from \eqref{faceini}, we obtain the result 
\be
F_{VV}(z;p,w)=\varphi(z;p)
\left(
\begin{array}{cccc}
1 & & & \\
 & X_{11}(z) & X_{12}(z) & \\
 & X_{21}(z) & X_{22}(z) & \\
 &  &  & 1\\
\end{array}
\right),
\en
where
\bea
&&\varphi(z;p)=
\frac{(pz;q^4,p)_\infty(pq^4z;q^4,p)_\infty}{(pq^2z;q^4,p)_\infty^2},
\lb{vphi}
\ena
and
\be
&&X_{11}(z)=\hg{wq^{2}}{q^{2}}{w}{p}{pq^{-2}z},
\\
&&X_{12}(z)=\frac{w(q-q^{-1})}{1-w}\hg{wq^{2}}{pq^{2}}{pw}{p}{pq^{-2}z},
\\
&&X_{21}(z)=\frac{pw^{-1}(q-q^{-1})}{1-pw^{-1}}z
\hg{pw^{-1}q^{2}}{pq^{2}}{p^2w^{-1}}{p}{pq^{-2}z},
\\
&&X_{22}(z)=\hg{pw^{-1}q^{2}}{q^{2}}{pw^{-1}}{p}{pq^{-2}z}.
\en
Here $\hg{q^a}{q^b}{q^c}{q}{z}$ denotes the basic hypergeometric series 
\be
\hg{q^a}{q^b}{q^c}{q}{z}=\sum_{n=0}^\infty
{(q^a;q)_n (q^b;q)_n \over (q;q)_n (q^c;q)_n} z^n.
\en

The image of the $R$ matrix  is determined from
\eqref{Rslth2} and the connection formula for the basic 
hypergeometric series
\begin{eqnarray*}
\hg{q^a}{q^b}{q^c}{q}{{1\over z}}
&=& 
{\Gamma_q(c)\Gamma_q(b-a) \Theta_q(q^{1-a}z)\over 
 \Gamma_q(b)\Gamma_q(c-a) \Theta_q(q z)}
\hg{q^a}{q^{a-c+1}}{q^{a-b+1}}{q}{q^{c-a-b+1}z}\\
&&+ 
{\Gamma_q(c)\Gamma_q(a-b) \Theta_q(q^{1-b}z)\over 
 \Gamma_q(a)\Gamma_q(c-b) \Theta_q(q z)}
\hg{q^b}{q^{b-c+1}}{q^{b-a+1}}{q}{q^{c-a-b+1}z},
\end{eqnarray*}
where
\begin{eqnarray*}
\Gamma_q(z)={(q;q)_\infty \over (q^z;q)_\infty }(1-q)^{1-z},\qquad
\Theta_q(z)= (z;q)_\infty (qz^{-1};q)_\infty (q;q)_\infty. 
\end{eqnarray*}
We find 
\be
&&R_{VV}(z;p,w)
=\rho(z;p)
\left(\begin{array}{cccc}
1 &&& \\
 & b(z;p,w) & c(z;p,w) & \\
 & \bar{c}(z;p,w)&  \bar{b}(z;p,w) & \\
 &&& 1 \end{array}
\right),
\en
with the coefficients given by 
\bea
&& 
\rho(z;p) = q^{-1/2}
\frac{(q^2 z;p,q^4)_\infty^2}{(z;p,q^4)_\infty(q^4 z;p,q^4)_\infty}
\frac{(p z^{-1};p,q^4)_\infty(pq^4 z^{-1};p,q^4)_\infty}
{(pq^2 z^{-1};p,q^4)_\infty^2},
\lb{rho}\\
&&b(z;p,w)= 
q{(pw^{-1}q^2;p)_\infty (pw^{-1}q^{-2};p)_\infty
\over (pw^{-1};p)_\infty^2}
{\Theta_p(z) \over \Theta_p(q^2z)},
\lb{bpw}\\
&&\bar{b}(z;p,w)= 
q{(wq^2;p)_\infty (wq^{-2};p)_\infty \over (w;p)_\infty^2}
{\Theta_p(z) \over \Theta_p(q^2z)},  
\lb{bbpw}\\
&& c(z;p,w)={\Theta_p(q^2)\over \Theta_p(w)}
{\Theta_p(wz)\over \Theta_p(q^2z)},
\lb{cpw}\\
&& \bar{c}(z;p,w)= z
{\Theta_p(q^2)\over \Theta_p(pw^{-1})}
{\Theta_p(pw^{-1}z)\over \Theta_p(q^2z)}.
\lb{cbpw}
\ena
As expected, these are (up to a gauge) the Boltzmann weights of the 
Andrews-Baxter-Forrester model \cite{ABF}. 

\subsection{The case $\Aqp{}(\slthbig)$} 

The case of $\Aqp{}(\slth)$ can be treated similarly. 
Let 
\bea
&&\ctR(\ze)=\Ad(\ze^\rho\otimes 1)(\cR),
\lb{Rslth3}
\\
&&E(\ze;p)=\Ad(\ze^\rho\otimes 1)(E(r)),
\lb{Eslth}
\\
&&\ctR(\ze;p)=\Ad(\ze^\rho\otimes 1)(\cR(r))
=\sigma\left(E(\ze^{-1};p)\right)\ctR(\ze)E(\ze;p)^{-1},
\lb{Rslth4}
\ena
where $p=q^{2r}$. 
In this case we have simply $q^T\ctR(0)=1$. 
Thus $E(\ze;p)$ is characterized by 
\be
&&E(p^{1/2}q^{c^{(1)}}\ze;p)=(\tau^{-1}\otimes\id)\left(E(\ze;p)\right)
\times q^{\widetilde{T}}\ctR(p^{1/2}q^{c^{(1)}}\ze),
\\
&&
E(0;p)=1.
\en

The calculation of the image in the two-dimensional representation 
can be done directly. 
Since $\pi\circ\tau=\Ad(\sigma^x)\circ\pi$ with 
$
\sigma^x=
\displaystyle\left(\begin{array}{cc} 0&1 \\ 1&0 \\ \end{array}\right)
$,
we have 
\be
&&(\pi\otimes\pi)(\tau^{2k}\otimes\id)\left(\ctR(p^k\ze)\right)
=q^{-1/2}\rho(p^{2k}\zeta^2)
\left(\begin{array}{cccc}
1 & & & \\
 & b_{2k}& c_{2k} & \\
 & c_{2k} & b_{2k} & \\
 &  &  & 1 \\
\end{array}\right),
\\
&&(\pi\otimes\pi)(\tau^{2k-1}\otimes\id)\left(\ctR(p^{k-\frac{1}{2}}\ze)\right)
=q^{-1/2}\rho(p^{2k-1}\zeta^2)
\left(\begin{array}{cccc}
b_{2k-1} & & & c_{2k-1} \\
 & 1&  & \\
 &  & 1 & \\
c_{2k-1} &  &  & b_{2k-1} \\
\end{array}\right),
\en
where $b_l,c_l$ are given in terms of \eqref{bR2} by
\be
b_l=b(p^l\ze^2),
\qquad 
c_l=p^{\frac{l}{2}}\ze\, c(p^l\ze^2).
\en
The infinite product can be readily calculated, yielding the result
\bea
&&
(\pi\otimes\pi)(E(\ze;p))
=
\varphi(\ze^2;p)
\left(\begin{array}{cccc}
a_E(\ze) & & & d_E(\ze)\\
 & b_E(\ze) & c_E(\ze) & \\
 & c_E(\ze) & b_E(\ze) & \\
d_E(\ze) &  &  & a_E(\ze) \\
\end{array}\right),
\lb{E1}
\ena
where $\varphi(z;p)$ is given by \eqref{vphi}, and 
\bea
&&
a_E(\ze)\pm d_E(\ze)
=\frac{(\mp p^{1/2}q\ze;p)_\infty}{(\mp p^{1/2}q^{-1}\ze;p)_\infty},
\lb{E2}\\
&&
b_E(\ze)\pm c_E(\ze)
=\frac{(\mp pq\ze;p)_\infty}{(\mp p q^{-1}\ze;p)_\infty}.
\lb{E3}
\ena
Finally the image of the $R$ matrix \eqref{Rslth4} is given by 
\bea
(\pi\otimes\pi)(\ctR(\ze;p))
=q^{-1/2}\rho(\zeta^2,p)
\left(\begin{array}{cccc}
a^+(\ze) & & & d^+(\ze)\\
 & b^+(\ze) & c^+(\ze) & \\
 & c^+(\ze) & b^+(\ze) & \\
d^+(\ze) &  &  & a^+(\ze) \\
\end{array}\right),
\lb{E4}
\ena
with 
\bea
&&a^+(\ze)\pm d^+(\ze)
=\frac{(\mp p^{1/2}q^{-1}\ze;p)_\infty}{(\mp p^{1/2}q\ze;p)_\infty}
\frac{(\mp p^{1/2}q\ze^{-1};p)_\infty}{(\mp p^{1/2}q^{-1}\ze^{-1};p)_\infty},
\lb{E5}\\
&&
b^+(\ze)\pm c^+(\ze)
=q\frac{1\pm q^{-1}\ze}{1\pm q\ze}
\frac{(\mp pq^{-1}\ze;p)_\infty}{(\mp p q\ze;p)_\infty}
\frac{(\mp pq\ze^{-1};p)_\infty}{(\mp p q^{-1}\ze^{-1};p)_\infty}.
\lb{E6}
\ena
This agrees with the 
$R$ matrix of the eight vertex model
(cf. eqs.(2.5)--(2.9) in \cite{hwm}, wherein 
the sign of $p^{1/2}$ is changed).
The results \eqref{E1}--\eqref{E6}
are due to Fr\o nsdal\cite{Fron,Fron2}.


\setcounter{section}{3}
\setcounter{equation}{0}
\section{Dynamical $RLL$-relations and vertex operators}

The $L$-operators and vertex operators for the elliptic algebras 
can be constructed from those of $U_q(\goth{g})$ 
by `dressing' the latter with the twistors. 
In this section, we examine various commutation relations among these 
operators. 
We shall mainly discuss 
the case of the face type algebra $\Bqla(\goth{g})$
where $\goth{g}$ is of affine type. 
We touch upon the vertex type algebras $\Aqp{}(\slnh)$ briefly at the end. 

\subsection{$RLL$-relation for $\Bqla(\gbig)$} 

Hereafter we write $U=U_q(\goth{g})$, $\B=\Bqla(\goth{g})$.  
By a representation of the quasi-Hopf algebra $\B$ we mean  
that of the underlying associative algebra $U$. 
Let $(\pi_V,V)$ be a finite dimensional module over $U$, 
and $(\pi_{V,z},V_z)$ be the evaluation representation associated with it
where $\pi_{V,z}=\pi_V\circ \Ad(z^d)$. 

\begin{dfn}
We define $L$-operators for $\B$ by 
\bea
&&L_V^\pm(z,\la)=\left(\pi_{V,z}\otimes \id \right)\cR^{'\pm}(\la),
\lb{Lpm}
\\
&&\cR^{'+}(\la)=q^{c\otimes d+d\otimes c}\cR(\la),
\lb{Rp}\\
&&\cR^{'-}(\la)=\cR^{(21)}(\la)^{-1}q^{-c\otimes d-d\otimes c}.
\lb{Rm}
\ena
Likewise we set
\be
&&
R^{\pm}_{VW}(z_1/z_2,\la)
=\left(\pi_{V,z_1}\otimes\pi_{W,z_2}\right)\cR^{'\pm}(\la).
\en
\end{dfn}

Setting further 
\bea
\cR^{'\pm}(z,\la)=\Ad(z^d\otimes 1)\cR^{'\pm}(\la), 
\lb{Rpmz}
\ena
we find from the dynamical YBE \eqref{dyYBE} that 
\be
&&
\cR^{'\pm(12)}(z_1/z_2,\la+h^{(3)})
\cR^{'\pm(13)}(q^{\mp c^{(2)}}z_1/z_3,\la)
\cR^{'\pm(23)}(z_2/z_3,\la+h^{(1)})
\\
&&\qquad \quad =
\cR^{'\pm(23)}(z_2/z_3,\la)
\cR^{'\pm(13)}(q^{\pm c^{(2)}}z_1/z_3,\la+h^{(2)})
\cR^{'\pm(12)}(z_1/z_2,\la),
\\
&&
\cR^{'+(12)}(q^{c^{(3)}}z_1/z_2,\la+h^{(3)})
\cR^{'+(13)}(z_1/z_3,\la)
\cR^{'-(23)}(z_2/z_3,\la+h^{(1)})
\\
&&\qquad \quad =
\cR^{'-(23)}(z_2/z_3,\la)
\cR^{'+(13)}(z_1/z_3,\la+h^{(2)})
\cR^{'+(12)}(q^{-c^{(3)}}z_1/z_2,\la).
\en
Applying $\pi_V\otimes\pi_W\otimes\id$, we obtain 
the dynamical $RLL$ relation. 
\begin{prop}\label{SOSRLL}
\bea
&&R^{\pm(12)}_{VW}(z_1/z_2,\la+ h)
L_V^{\pm(1)}(z_1,\la)L_W^{\pm(2)}(z_2,\la+ h^{(1)})
\nonumber\\
&&\quad =
L_W^{\pm(2)}(z_2,\la)L_V^{\pm(1)}(z_1,\la+ h^{(2)})
R^{\pm(12)}_{VW}(z_1/z_2,\la),
\lb{RLL1}
\\
&&R^{+(12)}_{VW}(q^cz_1/z_2,\la+ h)
L_V^{+(1)}(z_1,\la)L_W^{-(2)}(z_2,\la+ h^{(1)})
\nonumber\\
&&\quad =
L_W^{-(2)}(z_2,\la)L_V^{+(1)}(z_1,\la+ h^{(2)})
R^{+(12)}_{VW}(q^{-c}z_1/z_2,\la).
\lb{RLL2}
\ena
\end{prop}
Here the index ($1$) (resp. ($2$)) refers to $V$ (resp. $W$), and 
$h$, $c$ (without superfix) are elements of $\goth{h}\subset \B$. 
If we write 
\be
\la-\rho=rd+s'c+\bar{\la}-\bar{\rho} 
\qquad (r,s'\in \C, \bar{\la}\in\bar{\goth{h}}), 
\en
then 
\bea
\la+h^{(1)}=(r+h^{\vee}+c^{(1)})d+(s'+d^{(1)})c+
(\bar{\la}+\bar{h^{(1)}})
\lb{shift}
\ena
where $h^\vee$ is the dual Coxeter number. 
The parameter $r$ plays the role of the elliptic modulus.
Note that, in \eqref{RLL1}--\eqref{RLL2},  
$r$ also undergoes a shift 
depending on the central element $c$. 

Actually the two $L$-operators \eqref{Lpm} are not 
independent. 
\begin{prop}
We have 
\bea
L_V^+(pq^cz,\la)=q^{-2\bar{T}_{V,\bullet}}
(\Ad(\bar{X}_\la)\otimes\id)^{-1}L_V^-(z,\la),
\lb{LpLm}
\ena
where 
\be
&&\bar{T}_{V,\bullet}=\sum \pi(\bar{h}_j)\otimes \bar{h}^j,
\\
&&
\bar{X}_\la=\pi(q^{\sum \bar{h}_j\bar{h}^j+2(\bar{\la}-\bar{\rho})}).
\en
\end{prop}

\proof 
In the notation of \eqref{Rpmz}, we have
\be
&&
\cR^{'+}(z,\la)=
\sigma\left(F(z^{-1},\la)\right)\Bigl|_{z\mapsto q^{c^{(2)}-c^{(1)}}z}
\cR^{'+}(z)F(z,\la)^{-1},
\\
&&
\cR^{'-}(z,\la)=
\sigma\left(F(z^{-1},\la)\right)\cR^{'-}(z)
F(q^{c^{(2)}-c^{(1)}}z,\la)^{-1}.
\en
Now the difference equation \eqref{Fdiff} implies 
\be
&&q^{\bar{T}}\cR^{'+}(pz)F(pz,\la)^{-1}=
\left(\bar{\varphi}_{\la}\otimes\id\right)^{-1}
F(q^{-2c^{(1)}}z,\la),
\\
&&
\sigma\left(F(p^{-1}z^{-1},\la)\right)\Bigl|_{z\mapsto q^{-2c^{(1)}}z}
=
\left(\id\otimes\bar{\varphi}_\la\right)
\left(\sigma(F(z^{-1},\la))\cR^{'-}(z)q^{-\bar{T}}\right).
\en
Noting that 
\be
\Ad(q^{2\bar{T}})(\bar{\varphi}_\la\otimes\bar{\varphi}_\la)
\cR'=\cR',
\en
we find 
\be
\cR^{'+}(pq^{c^{(1)}+c^{(2)}}z,\la)
=q^{-2\bar{T}}\left(\bar{\varphi}_\la\otimes\id\right)^{-1}
\cR^{'-}(z,\la).
\en
Taking the image in $V$ we obtain the assertion.
\qed

\subsection{Vertex operators for $\Bqla(\gbig)$} 

Let $(\pi_{V,z},V_z)$ be as before, and let 
$V(\mu)$ be a highest weight module with highest weight $\mu$. 
Consider intertwiners of $U$-modules of the form 
\be
&&\Phi_V^{(\nu,\mu)}(z): V(\mu)\longrightarrow V(\nu)\otimes V_{z},
\\
&&\Psis_V^{(\nu,\mu)}(z): V_{z}\otimes  V(\mu)
\longrightarrow V(\nu).
\en
We call them vertex operators (VO's) of type I and type II, 
respectively.\footnote{
In this paper we treat VO's which have Fourier expansion in 
integral powers of $z$. In the notation of \cite{IIJMNT},
they are denoted with the symbol $~\widetilde{}~$; e.g. 
$\Phi_V^{(\nu,\mu)}(z)$ here is written as 
$\widetilde{\Phi}_\mu^{\nu,V}(z)$ in \cite{IIJMNT}.}
Define the corresponding VO's for $\B$ as follows \cite{Fron2}. 
\bea
&&\Phi_V^{(\nu,\mu)}(z,\la)=
(\id\otimes \pi_z)F(\la)\circ\Phi_V^{(\nu,\mu)}(z),
\label{twVO1} \\
&&\Psis_V^{(\nu,\mu)}(z,\la)=
\Psis_V^{(\nu,\mu)}(z)\circ (\pi_z\otimes \id)F(\la)^{-1}.
\label{twVO2}
\ena
When there is no fear of confusion, we often drop the sub(super)scripts 
$V$ or $(\nu,\mu)$. 
It is clear that \eqref{twVO1},\eqref{twVO2} 
satisfy the intertwining relations relative to the coproduct $\Delta_\la$
\eqref{facecopro}, 
\be
&&\De_\la(a) \Phi(z,\la)=\Phi(z,\la) a 
\qquad \forall a\in \B,
\\
&&a \Psis(z,\la)=\Psis(z,\la) \De_\la(a)
\qquad \forall a\in \B.
\en
These intertwining relations can be encapsulated to  
commutation relations with the $L$-operators. 

\begin{prop}\label{dynamicalLV}
The `dressed' VO's \eqref{twVO1},\eqref{twVO2} 
satisfy the following dynamical intertwining relations
(see the diagram below):
\bea
&&
\Phi_W(z_2,\la)L^+_V(z_1,\la)
=R^+_{VW}(q^cz_1/z_2,\la+h)L^+_V(z_1,\la)\Phi_W(z_2,\la+h^{(1)}),
\lb{dint1}\\
&&
\Phi_W(z_2,\la)L^-_V(z_1,\la)
=R^-_{VW}(z_1/z_2,\la+h)L^-_V(z_1,\la)\Phi_W(z_2,\la+h^{(1)}),
\lb{dint1m}
\\
&&L^+_V(z_1,\la)\Psi_W^*(z_2,\la+h^{(1)})
=\Psi_W^*(z_2,\la)L^+_V(z_1,\la+h^{(2)})R^+_{VW}(z_1/z_2,\la),
\lb{dint2}
\\
&&L^-_V(z_1,\la)\Psi_W^*(z_2,\la+h^{(1)})
=\Psi_W^*(z_2,\la)L^-_V(z_1,\la+h^{(2)})R^-_{VW}(q^cz_1/z_2,\la).
\lb{dint2m}
\ena
\end{prop}
%
\[
\begin{array}{ccc}
  V_{z_1}\otimes V(\mu)&\maprightu{1cm}{\Phi_W}\;\;
  V_{z_1}\otimes V(\nu)\otimes W_{z_2}\;\;\maprightu{1cm}{L^\pm_V}&
  V_{z_1}\otimes V(\nu)\otimes W_{z_2} \\
  \hspace{-0.4cm}\mapdownl{\sct L^\pm_V}&
  \makebox{\rule[-0.7cm]{0cm}{1.6cm}}&
  \mapdownr{\sct R^\pm_{VW}} \\
  V_{z_1}\otimes V(\mu)&\maprightd{5.5cm}{\Phi_W}&
  V_{z_1}\otimes V(\nu)\otimes W_{z_2}
\end{array}
\]

\medskip\medskip

\[
\begin{array}{ccc}
  V_{z_1}\otimes W_{z_2}\otimes V(\mu)\hspace{-0.7cm}&
  \maprightu{1cm}{R^\pm_{VW}}\;\;
  V_{z_1}\otimes W_{z_2}\otimes V(\mu)\;\;\maprightu{1cm}{L^\pm_V}&
  \hspace{-0.7cm}V_{z_1}\otimes W_{z_2}\otimes V(\mu) \\
  \hspace{0.1cm}\mapdownl{\sct \Psi^*_W}&
  \makebox{\rule[-0.7cm]{0cm}{1.6cm}}&
  \hspace{-1.1cm}\mapdownr{\sct \Psi^{*}_W} \\
  \hspace{0.4cm}V_{z_1}\otimes V(\nu)&\maprightd{7cm}{L^\pm_V}&
  \hspace{-0.8cm}V_{z_1}\otimes V(\nu)
\end{array}
\]
\medskip\medskip

\ignore{
\[
\begin{array}{ccccc}
V_{z_1}\otimes V(\mu)  & \stackrel{\Phi_W}{\longrightarrow} 
&V_{z_1}\otimes V(\nu)\otimes W_{z_2}&  \stackrel{L^\pm_V}{\longrightarrow}
&V_{z_1}\otimes V(\nu)\otimes W_{z_2}  \\
& &                &        &           \\
{\sct L^\pm_V}\downarrow \qquad   & & &       & \qquad\downarrow 
{\sct R^\pm_{VW}} \\
&&          &        &           \\
   V_{z_1}\otimes V(\mu)  & &\stackrel{\Phi_W}{\longrightarrow}
&   & V_{z_1}\otimes V(\nu)\otimes W_{z_2}
\end{array}
\]
\medskip

\[
\begin{array}{ccccc}
V_{z_1}\otimes W_{z_2}\otimes V(\mu)&
 \stackrel{R^\pm_{VW}}{\longrightarrow} &V_{z_1}\otimes W_{z_2}\otimes V(\mu)& 
 \stackrel{L^\pm_V}{\longrightarrow}
&V_{z_1}\otimes W_{z_2}\otimes V(\mu)  \\
& &                &        &           \\
{\sct \Psi^*_W}\downarrow\qquad    & & &       &\qquad \downarrow 
{\sct \Psi^{*}_W} \\
&&          &        &           \\
V_{z_1}\otimes V(\nu)  & &\stackrel{L^\pm_V}{\longrightarrow}
&   &V_{z_1}\otimes V(\nu)
\end{array}
\]
}

\proof
Let $\Delta$ be the original coproduct \eqref{copro} for $U$. 
The properties \eqref{qtri1}, \eqref{qtri2} are equivalently rewritten as 
\be
&&(\Delta\otimes\id)\cR(\la)=
F^{(12)}(\la+h^{(3)})^{-1}\cR^{(13)}(\la)\cR^{(23)}(\la+h^{(1)})
F^{(12)}(\la),
\\
&&(\id\otimes\Delta)\cR(\la)=
F^{(23)}(\la)^{-1}\cR^{(13)}(\la+h^{(2)})\cR^{(12)}(\la)
F^{(23)}(\la+h^{(1)}).
\en
{}From this it follows that 
\be
&&(\id\otimes\Delta)\cR^{'+}(z,\la)=
F^{(23)}(\la)^{-1}\cR^{'+(13)}(q^{c^{(2)}}z,\la+h^{(2)})
\cR^{'+(12)}(z,\la)F^{(23)}(\la+h^{(1)}),
\\
&&(\id\otimes\Delta)\cR^{'-}(z,\la)=
F^{(23)}(\la)^{-1}\cR^{'-(13)}(z,\la+h^{(2)})
\cR^{'-(12)}(q^{c^{(3)}}z,\la)F^{(23)}(\la+h^{(1)}).
\en
Using the intertwining relation 
\be
\Phi_W(z)a=\Delta(a)\Phi_W(z),
\en
we obtain 
\be
&&\phantom{=}\Phi_W(z_2,\la)L^+_V(z_1,\la)\\
&&=
\left(\pi_{V,z_1}\otimes\id\otimes\pi_{W,z_2}\right)
\left(F^{(23)}(\la)\right)\Phi_W(z_2)
\left(\pi_{V,z_1}\otimes \id\right)\left(\cR^{'+}(\la)\right)
\\
&&=
\left(\pi_{V,z_1}\otimes\id\otimes\pi_{W,z_2}\right)
\left(F^{(23)}(\la)(\id\otimes\Delta)\cR^{'+}(\la)\right)\Phi_W(z_2)
\\
&&=
\left(\pi_{V,z_1}\otimes\id\otimes\pi_{W,z_2}\right)
\left(\cR^{'+(13)}(q^{c^{(2)}},\la+h^{(2)})
\cR^{'+(12)}(\la)F^{(23)}(\la+h^{(1)})\right)
\Phi_W(z_2)
\\
&&=
R^+_{VW}(q^cz_1/z_2,\la+h)L^+_V(z_1,\la)\Phi_W(z_2,\la+h^{(1)}).
\en

The other cases are similar. \qed

{}From the theory of $q$KZ-equation\cite{FR}, 
we know the VO's for $U$ satisfy the 
commutation relations of the form 
\bea
&&\Rck_{VV}(z_1/z_2)\Phi_V^{(\nu,\mu)}(z_1)\Phi_V^{(\mu,\kappa)}(z_2)
=\sum_{\mu'}
\Phi_V^{(\nu,\mu')}(z_2)\Phi_V^{(\mu',\kappa)}(z_1)
W_{I}\BW{\kappa}{\mu}{\mu'}{\nu}{\frac{z_1}{z_2}},
\lb{VVt1}\\
&&\Psi_V^{*(\nu,\mu)}(z_1)\Psi_V^{*(\mu,\kappa)}(z_2)
{\Rck_{VV}(z_1/z_2)}^{-1}
=\sum_{\mu'}
W_{II}\BW{\kappa}{\mu}{\mu'}{\nu}{\frac{z_1}{z_2}}
\Psi_V^{*(\nu,\mu')}(z_2)\Psi_V^{*(\mu',\kappa)}(z_1),
\nonumber\\
&&\lb{VVt2}\\
&&\Phi_V^{(\nu,\mu)}(z_1)\Psi_V^{*(\mu,\kappa)}(z_2)
=\sum_{\mu'}
W_{I,II}\BW{\kappa}{\mu}{\mu'}{\nu}{\frac{z_1}{z_2}}
\Psi_V^{*(\nu,\mu')}(z_2)\Phi_V^{(\mu',\kappa)}(z_1).
\lb{VVt3}
\ena
Here 
\footnote{The $R$ matrix used here is the image of the universal $R$ matrix. 
It differs by a scalar factor from the one  used e.g. in \cite{JM,IIJMNT}
 in the commutation relations of VO's. }
\be
&&\Rck_{VV}(z)=PR_{VV}(z),
\qquad P(v\otimes v')=v'\otimes v,
\\
&&R_{VV}(z_1/z_2)=(\pi_{V,z_1}\otimes\pi_{V,z_2})\cR
\en
is the `trigonometric' $R$ matrix. 
In \eqref{VVt1}--\eqref{VVt3} we used a slightly abbreviated notation. 
For example, the left hand side of \eqref{VVt1} means the composition 
\be
V(\kappa)\maprightu{1cm}{\Phi(z_2)}V(\mu)\otimes V_{z_2}
\maprightu{2cm}{\Phi(z_1)\otimes \mbox{\footnotesize id}}
V(\nu)\otimes V_{z_1}\otimes V_{z_2}
\maprightu{2.5cm}{\mbox{\footnotesize id}\otimes\Rck(z_1/z_2)}
V(\nu)\otimes V_{z_1}\otimes V_{z_2}.
\en
Similarly \eqref{VVt2}, \eqref{VVt3} are maps 
\be
&&V_{z_2}\otimes V_{z_1}\otimes V(\kappa)\longrightarrow V(\nu),
\\
&&V_{z_2}\otimes V(\kappa) \longrightarrow V(\nu)\otimes V_{z_1},
\en
respectively.
For $U_q(\slth)$, 
the formulas for the $W$-factors in the simplest case can be found e.g. 
in \cite{IIJMNT}. 

The `dressed' VO's satisfy similar relations with appropriate dynamical 
shift. 
Setting $\Rck_{VV}(z,\la)=PR_{VV}(z,\la)$, we have 
\begin{prop}\label{dynamicalVV}
\bea
&&\Rck_{VV}(z_1/z_2,\la+h^{(1)})
\Phi_V^{(\nu,\mu)}(z_1,\la)\Phi_V^{(\mu,\kappa)}(z_2,\la)
\nn\\
&&\qquad\qquad 
=\sum_{\mu'}
\Phi_V^{(\nu,\mu')}(z_2,\la)\Phi_V^{(\mu',\kappa)}(z_1,\la)
W_{I}\BW{\kappa}{\mu}{\mu'}{\nu}{\frac{z_1}{z_2}},
\lb{VV11}\\
&&\Psi_V^{*(\nu,\mu)}(z_1,\la)\Psi_V^{*(\mu,\kappa)}(z_2,\la+h^{(1)})
\Rck_{VV}(z_1/z_2,\la)^{-1}
\nn\\
&&\qquad\qquad =\sum_{\mu'}
W_{II}\BW{\kappa}{\mu}{\mu'}{\nu}{\frac{z_1}{z_2}}
\Psi_V^{*(\nu,\mu')}(z_2,\la)\Psi_V^{*(\mu',\kappa)}(z_1,\la+h^{(1)}),
\lb{VV22}\\
&&\Phi_V^{(\nu,\mu)}(z_1,\la)\Psi_V^{*(\mu,\kappa)}(z_2,\la)
=\sum_{\mu'}
W_{I,II}\BW{\kappa}{\mu}{\mu'}{\nu}{\frac{z_1}{z_2}}
\Psi_V^{*(\nu,\mu')}(z_2,\la)\Phi_V^{(\mu',\kappa)}(z_1,\la+h^{(1)}).
\nn\\
&&\lb{VV12}
\ena
\end{prop}
Notice that the $W$-factors stay the same with the trigonometric case, 
and are not affected by a dynamical shift. 

\proof
Let us verify \eqref{VV22} as an example. 
We drop the suffix $V$. 
{}From the intertwining relation 
$a\Psis(z)=\Psis(z)\Delta(a)$, we have 
\bea
&&\mbox{LHS of \eqref{VV22}}
\nn\\
&&=
\Psi^{*(\nu,\mu)}(z_1)(\pi_z\otimes \id)(F(\la)^{-1})
\left(\id\otimes\Psi^{*(\mu,\kappa)}(z_2)\right)
\left(\pi_{z_1}\otimes\pi_{z_2}\otimes\id\right)(F^{(23)}(\la+h^{(1)})^{-1})
\nn\\
&&
=\Psi^{*(\nu,\mu)}(z_1)\left(\id\otimes\Psi^{*(\mu,\kappa)}(z_2)\right)
\nn\\
&&\qquad \times
\left(\pi_{z_1}\otimes\pi_{z_2}\otimes\id\right)
\left((\id\otimes \Delta)F(\la)^{-1}\cdot F^{(23)}(\la+h^{(1)})^{-1}\right).
\lb{aaa}
\ena
In the right hand side of \eqref{aaa},
the first factor equals 
\be
\sum_{\mu'}
W_{II}\BW{\kappa}{\mu}{\mu'}{\nu}{\frac{z_1}{z_2}}
\Psi^{*(\nu,\mu')}(z_2)\Psi^{*(\mu',\kappa)}(z_1)P^{(12)}
\left(\pi_{z_1}\otimes\pi_{z_2}\otimes\id\right)(\cR^{(12)}),
\en
while the second is 
\be
\left(\pi_{z_1}\otimes\pi_{z_2}\otimes\id\right)
\left((\Delta\otimes \id)F(\la)^{-1}F^{(12)}(\la)^{-1}\right)
\en
by the shifted cocycle condition. 
Since
\be
\cR^{(12)}(\Delta\otimes \id)F(\la)^{-1}F^{(12)}(\la)^{-1}
=
(\Delta'\otimes \id)F(\la)^{-1}F^{(21)}(\la)^{-1}\cR^{(12)}(\la),
\en
\eqref{aaa} becomes 
\be
&&
\sum_{\mu'}
W_{II}\BW{\kappa}{\mu}{\mu'}{\nu}{\frac{z_1}{z_2}}
\Psi^{*(\nu,\mu')}(z_2)\Psi^{*(\mu',\kappa)}(z_1)
\\
&&\times \left(\pi_{z_2}\otimes\pi_{z_1}\otimes\id\right)
\left((\Delta\otimes \id)F(\la)^{-1}F^{(12)}(\la)^{-1}\right)
\cdot P^{(12)}R(z_1/z_2,\la).
\en
Using again the shifted cocycle condition we arrive at the right 
hand side of \eqref{VV22}.
\qed

\subsection{The case of $\Aqp{}(\slnhbig)$} 

The case of vertex type algebras can be treated in a parallel way. 
Let $(\widetilde{\pi}_{V,\ze},V)$, 
$\widetilde{\pi}_{V,\ze}=\pi\circ\Ad(\ze^{\rho})$ 
stand for the evaluation module defined via the principal 
gradation operator $\rho$. 
In place of $d$ we use $\rho/n$ to define the $R$-matrix and 
$L$-operators as follows: 
\be
&&
\widetilde{L}_V^\pm(\ze,r)=\left(\widetilde{\pi}_{V,\ze}\otimes \id \right)
\widetilde{\cR}^{\,'\pm}(r),
\\
&&\widetilde{\cR}^{\,'+}(r)
=q^{\widetilde{T}}
\widetilde{\cR}(r),\\
&&\widetilde{\cR}^{\,'-}(r)=
\widetilde{\cR}^{(21)}(r)^{-1}
q^{-\widetilde{T}},
\\
&&
\widetilde{R}^{\pm}_{VW}(\ze_1/\ze_2,r)
=\left(\pi_{V,\ze_1}\otimes\pi_{W,\ze_2}\right)
\widetilde{\cR}^{\,'\pm}(r).
\en

Then the $RLL$ relations \eqref{RLL1}--\eqref{RLL2} remain 
valid if we replace the shift $\la+h$ by $r+c$ and read 
$q^c z$ as $q^{c/n} \ze$. 
Since $c$ is mapped to $0$ in $(\pi_{V,\ze},V)$, 
the relations somewhat simplify. 
The result reads as follows. 
\bea
&&\widetilde{R}^{\pm(12)}_{VW}(\ze_1/\ze_2,r+ c)
\widetilde{L}_V^{\pm(1)}(\ze_1,r)\widetilde{L}_W^{\pm(2)}(\ze_2,r)
\nonumber\\
&&\quad =
\widetilde{L}_W^{\pm(2)}(\ze_2,r)\widetilde{L}_V^{\pm(1)}(\ze_1,r)
\widetilde{R}^{\pm(12)}_{VW}(\ze_1/\ze_2,r),
\lb{RLL3}
\\
&&\widetilde{R}^{+(12)}_{VW}(q^{c/n}\ze_1/\ze_2,r+ c)
\widetilde{L}_V^{+(1)}(\ze_1,r)\widetilde{L}_W^{-(2)}(\ze_2,r)
\nonumber\\
&&\quad =
\widetilde{L}_W^{-(2)}(\ze_2,r)\widetilde{L}_V^{+(1)}(\ze_1,r)
\widetilde{R}^{+(12)}_{VW}(q^{-c/n}\ze_1/\ze_2,r).
\lb{RLL4}
\ena
By the same method as in the face type case, we find also
\be
\widetilde{\cR}^{\,'+}(p^{1/n}q^{(c^{(1)}+c^{(2)})/n}\ze,r)
=(\tau\otimes\id)^{-1}\widetilde{\cR}^{\,'-}(\ze,r).
\en
Taking the image in the vector representation 
$V=\C^n=\C v_1\oplus\cdots\oplus\C v_n$ and 
noting that $\tau$ is implemented by a conjugation
\be
\pi\circ\tau=\Ad(h)\circ \pi,
\qquad 
hv_j=v_{j+1\smod n},
\en
we obtain 
\bea
\widetilde{L}_V^+(p^{1/n}q^{c/n}\ze,r)
=(\Ad(h)\otimes\id)^{-1}\widetilde{L}_V^-(\ze,r).
\lb{Lrel}
\ena
The relations \eqref{RLL3},\eqref{RLL4} and \eqref{Lrel} 
first appeared (for $n=2$) in \cite{FIJKMY}.

Similarly we define the VO's by 
\be
&&\widetilde{\Phi}_V^{(\nu,\mu)}(\ze,r)=
(\id\otimes\pi_{V,\ze})(E(r))\circ\widetilde{\Phi}_V^{(\nu,\mu)}(\ze),
\\
&&\widetilde{\Psi}_V^{*(\nu,\mu)}(\ze,r)=
\widetilde{\Psi}_V^{*(\nu,\mu)}(\ze)\circ 
(\pi_{V,\ze}\otimes\id)(E(r))^{-1}.
\en
Here $\widetilde{\Phi}_V^{(\nu,\mu)}(\ze)$, 
$\widetilde{\Psi}_V^{*(\nu,\mu)}(\ze)$ 
are the VO's of $U_q(\slth)$ in the principal gradation.
The intertwining relations 
can be obtained from \eqref{dint1}--\eqref{dint2m}
by a simple replacement as explained above: 
\bea
&&
\widetilde{\Phi}_W(\ze_2,r)\widetilde{L}^+_V(\ze_1,r)
=
\widetilde{R}^+_{VW}(q^{c/n}\ze_1/\ze_2,r+c)\widetilde{L}^+_V(\ze_1,r)
\widetilde{\Phi}_W(\ze_2,r),
\lb{dint3}\\
&&
\widetilde{\Phi}_W(\ze_2,r)\widetilde{L}^-_V(\ze_1,r)
=\widetilde{R}^-_{VW}(\ze_1/\ze_2,r+c)\widetilde{L}^-_V(\ze_1,r)
\widetilde{\Phi}_W(\ze_2,r),
\lb{dint3m}
\\
&&\widetilde{L}^+_V(\ze_1,r)\widetilde{\Psi}_W^*(\ze_2,r)
=\widetilde{\Psi}_W^*(\ze_2,r)\widetilde{L}^+_V(\ze_1,r)
\widetilde{R}^+_{VW}(\ze_1/\ze_2,r),
\lb{dint4}
\\
&&\widetilde{L}^-_V(\ze_1,r)\widetilde{\Psi}_W^*(\ze_2,r)
=\widetilde{\Psi}_W^*(\ze_2,r)\widetilde{L}^-_V(\ze_1,r)
\widetilde{R}^-_{VW}(q^{c/n}\ze_1/\ze_2,r).
\lb{dint4m}
\ena
These formulas agree with those conjectured in \cite{FIJKMY,hwm}, 
if we identify $q^{2(r+c)}$ with $p$ there. 

The same can be done about the commutation relations of VO 
\eqref{VV11}--\eqref{VV12}. 
We do not repeat the formulas. 


\bigskip

{\it Acknowledgment.}\quad 
We thank
Hidetoshi Awata, 
Jintai Ding, 
Benjamin Enriquez, 
Boris Feigin, 
Ian Grojnowski, 
Koji Hasegawa,
Harunobu Kubo, 
Tetsuji Miwa ,
Takashi Takebe
and 
Jun Uchiyama
for discussions and interest. 

\appendix
\setcounter{equation}{0}
\section{Quasi-Hopf algebras}\lb{app:a}

We summarize here some basic notions concerning 
quasi-Hopf algebras \cite{QHA,QHA2}. 
Let $k$ be a commutative ring. 
\begin{dfn}
A quasi-bialgebra is a set $(A,\Delta,\varepsilon,\Phi)$ 
consisting of a unital associative $k$-algebra $A$, 
homomorphisms $\Delta: A \rightarrow A \otimes A$,  
$\varepsilon: A \rightarrow k$ 
and an invertible element $\Phi\in A \otimes A \otimes A$, 
satisfying the following axioms. 
\begin{eqnarray}
&& (\id \otimes \Delta)\Delta(a)=
\Phi (\Delta\otimes \id) \Delta(a) \Phi^{-1}
\quad \forall a\in A,\\
&&(\id \otimes \id \otimes \Delta) \Phi \cdot
  (\Delta \otimes \id \otimes \id) \Phi 
=
(1 \otimes \Phi) \cdot 
(\id \otimes  \Delta \otimes \id) \Phi \cdot
(\Phi \otimes 1),\\
&&(\varepsilon \otimes \id) \circ \Delta 
=\id=(\id \otimes \varepsilon) \circ \Delta,\\
&&(\id \otimes \varepsilon \otimes \id) \Phi=1.
\end{eqnarray}
\end{dfn}

\begin{dfn}
A quasi-Hopf algebra is a quasi-bialgebra $(A,\Delta,\varepsilon,\Phi)$ 
together with elements $\alpha,\beta \in A$ and 
an antiautomorphism $S$, satisfying the following conditions.
\begin{eqnarray}
\sum_i S(b_i) \alpha c_i=\ve(a) \alpha,
\qquad 
\sum_i b_i \beta S(c_i)=\ve(a) \beta, \label{qh1}
\end{eqnarray}
for $a\in A$, $\Delta(a)=\sum_i b_i \otimes c_i$, and 
\begin{eqnarray}
\sum_i X_i \beta S(Y_i)\alpha Z_i=1,\label{qh2}
\end{eqnarray}
where $\Phi=\sum_i X_i\otimes Y_i \otimes Z_i$.
\end{dfn}

\begin{dfn}
A quasi-triangular quasi-Hopf algebra is a set 
$(A,\Delta,\varepsilon,\Phi,R)$, 
where $(A,\Delta,\varepsilon,\Phi)$ is a quasi-Hopf algebra and 
$R\in A\otimes A$ is an invertible element such that 
\begin{eqnarray}
&& \Delta'(a)=R \Delta(a) R^{-1},
\lb{R1}\\
&& (\Delta\otimes \id) R=
\Phi^{(312)}R^{(13)} {\Phi^{(132)}}^{-1} R^{(23)}\Phi^{(123)},
\lb{R2}\\
&& (\id \otimes \Delta) R=
{\Phi^{(231)}}^{-1}R^{(13)} \Phi^{(213)} R^{(12)}{\Phi^{(123)}}^{-1}.
\lb{R3}
\end{eqnarray}
Here $\Delta'=\sigma\circ\Delta$ ($\sigma(a\otimes b)=b\otimes a$) 
is the opposite comultiplication.
\end{dfn}
In \eqref{R2}-\eqref{R3}, if $\Phi=\sum_i X_i\otimes Y_i \otimes Z_i$, 
then we write 
$\Phi^{(312)}=\sum_i Z_i\otimes X_i \otimes Y_i$, 
$\Phi^{(213)}=\sum_i Y_i\otimes X_i \otimes Z_i$, and so forth. 
Similar notations will be used throughout. 

The properties \eqref{R1}--\eqref{R3} imply in particular 
the Yang-Baxter type equation 
\begin{equation}
   R^{(12)} \Phi^{(312)} 
   R^{(13)} {\Phi^{(132)}}^{-1}   
   R^{(23)} \Phi^{(123)} 
=
  \Phi^{(321)} R^{(23)} 
  {\Phi^{(231)}}^{-1} R^{(13)}  
  \Phi^{(213)} R^{(12)}.
\lb{RPhi}
\end{equation}

There is an important operation called {\it twist}, 
which associates a new quasi-bialgebra with a given one. 
Let $(A,\Delta,\varepsilon,\Phi)$ be a quasi-bialgebra, 
and let $F\in A\otimes A$ be an invertible element such that  
$(\id\otimes \ve)F=1=(\ve\otimes \id)F$.
Set
\begin{eqnarray}
&&\tilde{\Delta}(a)=F \Delta(a) {F}^{-1}
\qquad (\forall a\in A), 
\label{qbt1} \\
&& \tilde{\Phi}=\left (F^{(23)} (\id \otimes \Delta) F \right) \,\Phi\, 
\left (F^{(12)} (\Delta \otimes \id) F \right)^{-1}.
\label{qbt2} 
\end{eqnarray}
Then $(A,\tilde{\Delta},\varepsilon,\tilde{\Phi})$ is also a quasi-bialgebra. 
We refer to the element $F$ as {\it twistor}. 
If in addition $(A,\Delta,\varepsilon,\Phi)$ is a quasi-Hopf algebra, 
with $\alpha,\beta,S$ satisfying (\ref{qh1}) and (\ref{qh2}), 
then $(A,\tilde{\Delta},\varepsilon,\tilde{\Phi})$ 
defined by  (\ref{qbt1}) and (\ref{qbt2}) together with 
\be
\tilde{S}=S, 
\qquad
\tilde{\alpha}=\sum_i S(d_i) \alpha e_i,
\qquad 
\tilde{\beta}=\sum_i f_i\beta S(g_i), 
\en
is also a quasi-Hopf algebra.  
Here we have set $\sum_i d_i\otimes e_i=F^{-1}$ and $\sum_i f_i\otimes g_i=F$.
Finally, a twist of a 
quasi-triangular quasi-Hopf algebra is again quasi-triangular, 
with the choice of new $R$ given by 
\begin{eqnarray}
\tilde{R}=F^{(21)} R {F^{(12)}}^{-1}. 
\end{eqnarray}

An important special case is a twist of a quasi-triangular {\it Hopf} algebra
$(A,\Delta,\varepsilon,R)$ 
(i.e. a quasi-triangular quasi-Hopf algebra with $\Phi=1$)
by a shifted cocycle. 
Let $H$ be an abelian subalgebra of $A$, with the product written additively.

\begin{dfn} A twistor $F(\la)$ depending on $\la\in H$ is a shifted cocycle 
if it satisfies the relation 
\bea
F^{(12)}(\la)\left(\Delta\otimes\id\right)F(\la)
=F^{(23)}(\la+h^{(1)})\left(\id\otimes\Delta\right)F(\la)
\lb{cocyA} 
\ena
for some $h\in H$. 
\end{dfn}

Let $(A,\Delta_\la,\varepsilon,\Phi(\la),R(\la))$ be the 
quasi-triangular quasi-Hopf algebra obtained by a twist by $F(\la)$. 
The shifted cocycle condition \eqref{cocyA} simplifies 
the properties of $\Phi(\la)$ and $R(\la)$ as follows.

\begin{prop} We have 
\bea
&&\Phi(\la)=F^{(23)}(\la)F^{(23)}(\la+h^{(1)})^{-1},
\lb{qtri0}\\
&&
(\Delta_{\la}\otimes\id) R(\la)
= \Phi^{(312)}(\la) R^{(13)}(\la)R^{(23)}(\la+ h^{(1)}),
\lb{qtri1}\\
&&
(\id \otimes \Delta_{\la}) R(\la)
=
R^{(13)}(\la+ h^{(2)})R^{(12)}(\la)  {\Phi^{(123)}}(\la)^{-1}.
\lb{qtri2}
\ena
As a corollary the dynamical Yang-Baxter relation holds:
\begin{equation}
R^{(12)}(\la+ h^{(3)} )R^{(13)}(\la)R^{(23)}(\la+h^{(1)})
=R^{(23)}(\la)R^{(13)}(\la+h^{(2)})R^{(12)}(\la).
\label{dyYBE}
\end{equation}
\end{prop}



\begin{thebibliography}{10}

\bibitem{FIJKMY}
O.~Foda, K.~Iohara, M.~Jimbo, R.~Kedem, T.~Miwa, and H.~Yan.
\newblock An elliptic quantum algebra for $\widehat{sl}_2$.
\newblock {\em Lett. Math. Phys.}, 32:259--268, 1994.

\bibitem{Fel95}
G.~Felder.
\newblock Elliptic quantum groups.
\newblock {\em Proc. {ICMP Paris} 1994}, pages 211--218, 1995.

\bibitem{Fron}
C.~Fr\o nsdal.
\newblock Generalization and exact deformations of quantum groups.
\newblock {\em Publ.RIMS, Kyoto Univ.}, 33:91--149, 1997.

\bibitem{Fron1}
C.~Fr\o nsdal.
\newblock Quasi-{Hopf} deformation of quantum groups.
\newblock {\em Lett. Math. Phys.}, 40:117--134, 1997.

\bibitem{Bax72}
R.~J. Baxter.
\newblock Partition function of the eight-vertex lattice model.
\newblock {\em Ann. of Phys.}, 70:193--228, 1972.

\bibitem{Bel}
A.~Belavin.
\newblock Dynamical symmetry of integrable quantum systems.
\newblock {\em Nucl. Phys.}, B 180 [FS2]:189--200, 1981.

\bibitem{Skl82}
E.~K. Sklyanin.
\newblock Some algebraic structure connected with the {Yang-Baxter} equation.
\newblock {\em Funct. Anal. Appl.}, 16:263--270, 1982.

\bibitem{Cher85}
I.~V. Cherednik.
\newblock Some finite-dimensional representations of generalized {Sklyanin}
  algebras.
\newblock {\em Funct. Anal. Appl.}, 19:77--79, 1985.

\bibitem{FeOd95}
B.~L. Feigin and A.~V. Odesskii.
\newblock Vector bundles on elliptic curve and {Sklyanin} algebras, 1995.
\newblock RIMS-1032, q-alg/9509021.

\bibitem{ABF}
G.~E. Andrews, R.~J. Baxter, and P.~J. Forrester.
\newblock Eight-vertex {SOS} model and generalized {Rogers-Ramanujan}-type
  identities.
\newblock {\em J. Stat. Phys.}, 35:193--266, 1984.

\bibitem{fusion}
E.~Date, M.~Jimbo, T.~Miwa, and M.~Okado.
\newblock Fusion of the eight vertex {SOS} model.
\newblock {\em Lett. Math. Phys.}, 12:209--215, 1986.

\bibitem{JMO1}
M.~Jimbo, T.~Miwa, and M.~Okado.
\newblock Solvable lattice models whose states are dominant integral weights of
  {$A^{(1)}_{n-1}$}.
\newblock {\em Lett. Math. Phys.}, 14:123--131, 1987.

\bibitem{JMO3}
M.~Jimbo, T.~Miwa, and M.~Okado.
\newblock Solvable lattice models related to the vector representation of
  classical simple {Lie} algebras.
\newblock {\em Commun. Math. Phys.}, 116:507--525, 1988.

\bibitem{BBB}
O.~Babelon, D.~Bernard, and E.~Billey.
\newblock A quasi-{Hopf} algebra interpretation of quantum {$3j$}- and
  {$6j$-}symbols and difference equations.
\newblock {\em Phys. Lett. B}, 375:89--97, 1996.

\bibitem{QHA}
V.~G. Drinfeld.
\newblock Quasi-{Hopf} algebras.
\newblock {\em Leningrad Math. J.}, 1:1419--1457, 1990.

\bibitem{EF}
B.~Enriquez and G.~Felder.
\newblock Elliptic quantum groups {$E_{\tau,\eta}(\slth)$} and quasi-{Hopf}
  algebras, 1997.
\newblock q-alg/9703018.

\bibitem{Konno}   
H.~Konno.   
\newblock An elliptic algebra $U_{q,p}\bigl(\slth\bigr)$
and the fusion {RSOS} model.
\newblock {\em Commun. Math. Phys.}, 195:373-403, 1998.

\bibitem{Fa82}
L.~D. Faddeev.
\newblock Integrable models in {$(1+1)$}-dimensional quantum field theory.
\newblock {\em Les Houches Lectures}, {XXXIX}:561--608, 1982.

\bibitem{ReFa83}
N.~Reshetikhin and L.~Faddeev.
\newblock Hamiltonian structures for integrable models of field theory.
\newblock {\em Theoret. Math. Phys.}, 56:847--862, 1983.

\bibitem{Ta85}
L.~A. Takhtadzhan.
\newblock Solutions of the triangle equations with {${\bf Z}_n\times {\bf
  Z}_n$}-symmetry and matrix analogs of the {Weierstrass} zeta- and sigma-
  functions.
\newblock {\em J. Soviet Math.}, 31:3432--3444, 1985.

\bibitem{Dri86}
V.~G. Drinfeld.
\newblock Quantum groups.
\newblock {\em Proc. ICM, Am. Math. Soc., Berkeley, CA}, pages 798--820, 1986.

\bibitem{Tani}
T.~Tanisaki.
\newblock Killing forms, {Harish-Chandra} isomorphisms, and universal
  {$R$}-matrices for quantum algebras.
\newblock {\em Int. J. Mod. Phys.}, A7 supplement 1B:941--961, 1992.

\bibitem{IIJMNT}
M.~Idzumi, K.~Iohara, M.~Jimbo, T.~Miwa, T.~Nakashima, and T.~Tokihiro.
\newblock Quantum affine symmetry in vertex models.
\newblock {\em Int. J. Mod. Phys.}, A8:1479--1511, 1993.

\bibitem{hwm}
O.~Foda, K.~Iohara, M.~Jimbo, R.~Kedem, T.~Miwa, and H.~Yan.
\newblock Notes on highest weight modules of the elliptic algebra {${\cal
  A}_{q,p}(\widehat{sl}_2)$}.
\newblock {\em Prog. Theoret. Phys., Supplement}, 118:1--34, 1995.

\bibitem{Fron2}
C.~Fr\o nsdal and A.~Galindo.
\newblock {$8$}-vertex correlation functions and twist covariance of {$q$-KZ}
  equation.
\newblock q-alg/9709028.

\bibitem{FR}
I.~B.~Frenkel and N.~Yu~Reshetikhin,
\newblock Quantum affine algebras and holonomic difference equations.
\newblock {\em Commun. Math. Phys.}, 146:1--60, 1992.

\bibitem{JM}
M.~Jimbo and T.~Miwa.
\newblock {\em Algebraic Analysis of Solvable Lattice Models}.
\newblock CBMS Regional Conference Series in Mathematics vol. 85, AMS, 1994.

\bibitem{QHA2}
V.~G. Drinfeld.
\newblock On quasitriangular quasi-{Hopf} algebras and a group closely
  connected with {Gal$(\overline{{\bf Q}}/{\bf Q})$}.
\newblock {\em Leningrad Math. J.}, 2:829--860, 1991.

\end{thebibliography}

\end{document}